\documentclass[12pt]{article} \usepackage{amssymb}

\newtheorem{thm}{Theorem}[section]

\newtheorem{prop}[thm]{Proposition}
\newtheorem{cor}[thm]{Corollary}
\newtheorem{rema}[thm]{Remark}
\newtheorem{app}[thm]{Application}
\newtheorem{prob}[thm]{Problem}
\newcommand{\halmos}{\rule{1ex}{1.4ex}}
\newcommand{\pfbox}{\hspace*{\fill}\mbox{$\halmos$}}

\begin{document}

\begin{center}
\begin{Large}
{\bf Factorization of formal exponentials and uniformization}
\end{Large}

\medskip

Katrina Barron\footnote{Supported in part by an NSF Mathematical
Sciences Postdoctoral Research Fellowship and by a University of
California President's Postdoctoral Fellowship}, 
\\ Department of Mathematics, University of California, Santa Cruz, 
CA 95064\\

Yi-Zhi Huang\footnote{Supported in part by NSF grant DMS-9622961} 
\\ Institut des Hautes \'{E}tudes Scientifiques, F-91440
Bures-sur-Yvette, France and Department of Mathematics, Rutgers
University, Piscataway, NJ 08854\\

and \\

James Lepowsky\footnote{Supported in part by NSF grants
DMS-9401851 and DMS-9701150}\\
Department of Mathematics, Rutgers University, Piscataway, NJ 08854
\end{center}

\hspace{1.5 cm}
\begin{abstract}
Let $\mathfrak{g}$ be a Lie algebra over a field of characteristic
zero equipped with a vector space decomposition $\mathfrak{g} =
\mathfrak{g}^- \oplus \mathfrak{g}^+$, and let $s$ and $t$ be
commuting formal variables commuting with $\mathfrak{g}$.  We prove
that the map $C: s\mathfrak{g}^-[[s,t]] \times t\mathfrak{g}^+ [[s,t]]
\longrightarrow s\mathfrak{g}^-[[s,t]] \oplus t\mathfrak{g}^+ [[s,t]]$
defined by the Campbell-Baker-Hausdorff formula and given by $e^{sg^-}
e^{tg^+} = e^{C(sg^-,tg^+)}$ for $g^\pm \in \mathfrak{g}^\pm [[s,t]]$
is a bijection, as is well known when $\mathfrak{g}$ is
finite-dimensional over $\mathbb{R}$ or $\mathbb{C}$, by geometry.  It
follows that there exist unique $\Psi^\pm \in \mathfrak{g}^\pm[[s,t]]$
such that $e^{tg^+} e^{sg^-} = e^{s\Psi^-} e^{t\Psi^+}$ (also well
known in the finite-dimensional geometric setting).  We apply this to
a Lie algebra $\mathfrak{g}$ consisting of certain formal infinite
series with coefficients in a $\mathbb{Z}$-graded Lie algebra
$\mathfrak{p}$, for instance, an affine Lie algebra, the Virasoro
algebra or a Grassmann envelope of the $N = 1$ Neveu-Schwarz
superalgebra.  For $\mathfrak{p}$ the Virasoro algebra, the result was
first proved by Huang as a step in the construction of a geometric
formulation of the notion of vertex operator algebra, and for
$\mathfrak{p}$ a Grassmann envelope of the Neveu-Schwarz
superalgebra, it was first proved by Barron as a corresponding step in
the construction of a supergeometric formulation of the notion of
vertex operator superalgebra.  In the special case of the Virasoro
(resp., $N=1$ Neveu-Schwarz) algebra with zero central charge the
result gives the precise expansion of the uniformizing function for a
sphere (resp., supersphere) with tubes resulting {}from the sewing of
two spheres (resp., superspheres) with tubes in two-dimensional
genus-zero holomorphic conformal (resp., $N = 1$ superconformal) field
theory.  The general result places such uniformization problems into a
broad formal algebraic context.
\end{abstract}

\section{Introduction}

Recall that for a Lie algebra $\mathfrak{g}$ over a field of
characteristic zero and $g_1, g_2 \in \mathfrak{g}$, the classical
Campbell-Baker-Hausdorff formula (cf. \cite{R}) gives a formal Lie
series $C(g_1,g_2)$ in $g_1$ and $g_2$ such that $e^{g_1} e^{g_2} =
e^{C(g_1,g_2)}$.  In this paper we prove the following factorization
theorem (see Theorem \ref{main theorem} below): Let $\mathfrak{g}$ be
equipped with a vector space decomposition $\mathfrak{g} =
\mathfrak{g}^- \oplus \mathfrak{g}^+$. Then the map
\begin{equation}\label{C in intro}
C: s\mathfrak{g}^-[[s,t]] \times t\mathfrak{g}^+[[s,t]]
\longrightarrow s\mathfrak{g}^-[[s,t]] \oplus t\mathfrak{g}^+[[s,t]]
\end{equation}
defined by the Campbell-Baker-Hausdorff formula and given by 
\begin{equation}\label{main equation in intro}
e^{sg^-} e^{tg^+} = e^{C(sg^-,tg^+)}
\end{equation} 
for $g^\pm \in \mathfrak{g}^\pm [[s,t]]$ is a bijection; here $s$ and 
$t$ are commuting formal variables commuting with $\mathfrak{g}$.  
The map $C$ is just the Campbell-Baker-Hausdorff formula inside the 
Lie algebra $\mathfrak{g}[[s,t]]$, the coefficients $s$ and $t$ in 
(\ref{main equation in intro}) making this possibly infinite Lie series 
well-defined.  The content of our theorem is that $C^{-1}$ exists and 
is a factorization of formal exponentials with respect to this vector 
space decomposition.  (Note that in the domain of the map $C$ in 
(\ref{C in intro}), we have the cartesian product of the two spaces, 
while in the codomain we of course have the same set but viewed as the 
vector space direct sum.)

In the case that $\mathfrak{g}$ is a finite-dimensional Lie algebra
over $\mathbb{R}$ or $\mathbb{C}$, this result is well known and is
proved using the geometry of a corresponding Lie group (see
e.g. \cite{V}). But in the case of infinite-dimensional Lie algebras,
since the theory of the corresponding group-like structures and of the
correspondence between Lie algebras and these group-like structures is
not developed, the argument proving in the finite-dimensional Lie
algebra case that for given $h^\pm \in \mathfrak{g}^\pm [[s,t]]$, one 
can factor $e^{sh^- + th^+}$ uniquely as $e^{sg^-} e^{tg^+}$ for some 
$g^\pm \in \mathfrak{g}^\pm [[s,t]]$ cannot be generalized directly.  
Even more to the point, in the case of the Virasoro algebra and the 
$N=1$ Neveu-Schwarz (superconformal) algebra, special cases of this 
result were indeed needed and proved in the very study of the 
correspondence between the infinite-dimensional Lie algebras and the 
group-like structures (see \cite{H1}-\cite{H3} and
\cite{B1}-\cite{B2}).  To our knowledge, the present general formal 
factorization result in the infinite-dimensional case has not 
previously appeared in the literature.  Regardless, in the present 
work we are mainly concerned with applications of the result (in the 
infinite-dimensional case).  

We also give precise information about the form of the elements 
$sg^-$ and $tg^+$ in terms of $sh^-$ and $th^+$ by defining 
universal formal series that we call the
``canonical factorization series'' $F^\pm$, which are certain 
formal infinite linear combinations of ``words''
in $sh^-$, $th^+$ and the canonical projections $\pi^\pm : 
\mathfrak{g}[[s,t]] \rightarrow \mathfrak{g}^\pm [[s,t]]$, and
showing that $sg^- = F^-(sh^-,th^+;\pi^\pm)$ and $tg^+ = 
F^+(sh^-,th^+;\pi^\pm)$.

As a corollary of the factorization theorem stated above, we use
(\ref{C in intro}) and (\ref{main equation in intro}) to construct
(see Corollary \ref{switching corollary} below) a unique bijection:
\begin{equation}
\Psi : t\mathfrak{g}^+[[s,t]] \times s\mathfrak{g}^-[[s,t]]
\longrightarrow s\mathfrak{g}^-[[s,t]] \times t\mathfrak{g}^+[[s,t]]
\end{equation}
such that for $g^\pm \in \mathfrak{g}^\pm[[s,t]]$,
\begin{equation}
 e^{tg^+} e^{sg^-} = e^{s\Psi^-} e^{t\Psi^+}  
\end{equation}
with $s\Psi^- = \pi^- \circ \Psi (tg^+,sg^-)$ and $t\Psi^+ = \pi^+
\circ \Psi (tg^+,sg^-)$.  We call this result ``formal algebraic 
uniformization,'' for reasons about to be explained.  It follows 
{}from the information contained in the canonical factorization
series arising {}from the factorization theorem, and the
Campbell-Baker-Hausdorff theorem, that $s\Psi^-$ and
$t\Psi^+$ are equal to certain universal formal series of words
in $sg^-$, $tg^+$ and the projections $\pi^\pm$.  We call these 
canonical series the ``formal algebraic uniformization series.''

In our applications, which we are about to describe, we use the
existence of the canonical factorization series and of the formal
algebraic uniformization series as steps in the proof.  These formal
series of words can be thought of as analogues, in some sense, of the
Campbell-Baker-Hausdorff series.
 
For $\mathfrak{p}$ a $\mathbb{Z}$-graded Lie algebra, 
we apply the formal algebraic uniformization (and also factorization)
results
to a Lie subalgebra $\mathfrak{g}$ of the Lie algebra
consisting of certain formal
infinite series with coefficients in $\mathfrak{p}$ (Corollary
\ref{application corollary} below).  In \cite{H1}, Huang proved
Corollary \ref{application corollary} in the case where $\mathfrak{p}$
is the Virasoro algebra (see Application \ref{virasoro example} below)
by first establishing the result in a certain representation of the
Virasoro algebra, namely, the standard representation given by $L_n
\mapsto - x^{n + 1} \frac{d}{dx} \in \mathrm{End} (\mathbb{C}[x,
x^{-1}])$ and $c = 0$, and then lifting to representations with
general $c$, and extending to a modification of the universal
enveloping algebra.  In \cite{B1}, Barron used a similar approach to
prove Corollary \ref{application corollary} in the case where
$\mathfrak{p}$ is a Grassmann envelope of the $N = 1$ Neveu-Schwarz
superalgebra (see Application \ref{N = 1} below).  These two cases are
fundamental to the sewing operations in conformal and $N = 1$
superconformal field theory, respectively (see \cite{H1}-\cite{H3},
\cite{B1}-\cite{B2}).  In the case of the Virasoro (resp., $N=1$
Neveu-Schwarz) algebra with zero central charge the result gives the
precise expansion of the uniformizing function for a sphere (resp.,
supersphere) with tubes resulting {}from the sewing of two spheres
(resp., superspheres) with tubes in two-dimensional genus-zero
holomorphic conformal (resp., $N = 1$ superconformal) field
theory. This paper gives a unified proof of these two results and
further shows that these are special cases of a much more general
result, namely, Corollary \ref{application corollary}. This corollary
can in addition be applied to obtain the corresponding results for a
Grassmann envelope of the $N > 1$ Neveu-Schwarz algebras of
supersymmetries in $N > 1$ superconformal field theories, and also to
such structures as affine Lie algebras (see Application \ref{affine}
below) and superalgebras.  Since Corollary \ref{application corollary}
is in turn a special case of Corollary \ref{switching corollary}, we
can perhaps view Corollary \ref{switching corollary} as a generalized
and canonical formal algebraic version of such uniformization.

To prove the results of the present paper, one might have hoped that
the method used in \cite{H1} and \cite{B1} to obtain the results in
special representations could be generalized directly to the universal
enveloping algebra arising in the formulations of our main results. 
However, the direct generalization of that method to the universal 
enveloping algebra does not work because the method in \cite{H1} and 
\cite{B1} uses properties of the special representations that 
universal enveloping algebras do not have. In the present paper, 
instead of working in the universal enveloping algebra, we work 
directly in the Lie algebra (see Remark \ref{universal vs lie} below). 
In particular, we reformulate the desired results as equations in the 
Lie algebra, and the method in \cite{H1} and \cite{B1} used for 
special representations can now be applied in the Lie algebra to solve 
these equations and thus obtain the results (see Remark \ref{method of 
proof} below).  In fact a crucial observation in this paper is that 
even though there is no setting involving universal enveloping algebras 
to which the method of proof used in \cite{H1} and \cite{B1} can be 
applied, one can in fact still find a very general setting to which the 
method can be applied, a setting very different {}from the ones in 
\cite{H1} and \cite{B1}.

This paper is organized as follows.  In Section 2, we give a review of
the Campbell-Baker-Hausdorff formula, including some notation that
will be useful later.  In Section 3, we prove the main theorem on the
factorization of formal exponentials, we establish a number of related
results, and we use the main theorem to prove
the corollary giving formal algebraic uniformization.  
We also introduce the canonical factorization series and the formal
algebraic uniformization series.  In Section 4,
we apply these results to Lie algebras consisting of certain formal
infinite series with coefficients in a $\mathbb{Z}$-graded Lie
algebra.  We then give several applications of the result for affine
Lie algebras, the Virasoro algebra and Grassmann envelopes of the $N =
1$ Neveu-Schwarz superalgebra.  For the latter two applications, we
discuss the importance of these results to conformal and
superconformal field theory, and we point out that the result also
applies to Grassmann envelopes of the $N>1$ Neveu-Schwarz algebras.

\section{The Campbell-Baker-Hausdorff formula}

We begin with some review of the Campbell-Baker-Hausdorff formula,
following \cite{R}.  Let $\mathbb{F}$ be a field of characteristic
zero, and let $a$ and $b$ be two formal noncommuting symbols.  Let
$\mathbb{F} \langle a,b \rangle$ be the $\mathbb{F}$-algebra of formal
linear combinations of words in $a$ and $b$, i.e., noncommutative
polynomials in $a$ and $b$ over $\mathbb{F}$.  A {\em Lie polynomial}
in $a$ and $b$ is an element of the smallest $\mathbb{F}$-subspace of
$\mathbb{F} \langle a,b \rangle$ containing $a$ and $b$ and closed
under Lie brackets.  For an element $S$ of the formal completion
$\mathbb{F} \langle \langle a,b \rangle \rangle$ of $\mathbb{F}
\langle a,b \rangle$, write
\[S = \sum_{n \in \mathbb{N}} S_n ,\]
where each $S_n$ is homogeneous of total degree $n$ in $a$ and $b$.
Then $S$ is called a {\em Lie series} if each $S_n$ is a Lie
polynomial.  For any formal series $S \in \mathbb{F} \langle \langle
a,b \rangle \rangle$ with no constant term, its formal exponential 
\[e^S = \sum_{n \in \mathbb{N}} \frac{S^n}{n!} \]
is well defined.

The Campbell-Baker-Hausdorff theorem asserts the existence of a
(unique) Lie series $C(a,b) \in \mathbb{Q} \langle \langle a,b \rangle
\rangle$ such that
\begin{equation}\label{CH}
e^a e^b = e^{C(a,b)} . 
\end{equation}

We now recall the precise formula for $C(a,b)$.  Even though we will
not need the main information contained in it, it is nice to see the
role that it plays in our proof of the bijectivity of $C$ in Theorem
\ref{main theorem} below.  For $S \in \mathbb{Q} \langle \langle a,b
\rangle \rangle$, let $S \frac{\partial}{\partial b}$ denote the
derivation of $\mathbb{Q} \langle \langle a,b \rangle \rangle$ that
maps $a$ to $0$ and $b$ to $S$. The series $C(a,b)$ is given by
\begin{equation}\label{real CH formula}
 C(a,b) = \exp \left( H_1 \frac{\partial}{\partial b} \right) (b), 
\end{equation}
where 
\[H_1 = \left(\frac{\mathrm{ad} \; b}{e^{\mathrm{ad} \; b} - 1}
\right)(a) ,\] 
and $H_1$ is the part of $C(a,b)$ that is homogeneous of degree one
with respect to $a$.  In particular, if we write
\begin{equation}\label{C decomp 1}
C(a,b) = \sum_{j \in \mathbb{N}} H_j 
\end{equation}
where $H_j$ is the part of $C(a,b)$ that is homogeneous of degree $j$ 
with respect to $a$, then
\[ H_j = \frac{1}{j!} \left( H_1 \frac{\partial}{\partial b} \right)^j
(b) \] 
(cf. \cite{R}).  Furthermore, writing
\begin{equation}\label{C decomp 2}
H_j = \sum_{k \in \mathbb{N}} H_{j,k},  
\end{equation} 
where $H_{j,k} \in \mathbb{Q} \langle a,b \rangle$ is homogeneous of
degree $j$ in $a$ and degree $k$ in $b$, we note that
\begin{equation}\label{initial conditions for H}
H_{j,0} = \left\{ \begin{array}{ll} 
               a  & \mbox{if $j = 1$}\\
               0  & \mbox{otherwise}
\end{array}  \right., \qquad
H_{0,k} = \left\{ \begin{array}{ll} 
               b  & \mbox{if $k = 1$}\\
               0  & \mbox{otherwise}
\end{array}  \right. 
\end{equation}
and 
\begin{equation}\label{second initial conditions for H} 
H_{1,1} = \frac{1}{2} [a,b] .
\end{equation}

\section{Factorization of formal exponentials and formal algebraic 
uniformization}

In this section we state and prove the two main results of the paper:
factorization of formal exponentials (Theorem \ref{main theorem}) and
formal algebraic uniformization (Corollary \ref{switching corollary}).
We also give precise information about the form of the resulting
elements (Theorem \ref{generalwordsthm} and Corollary
\ref{switchingwordscorollary}).

We work over our field $\mathbb{F}$ of characteristic zero.  We fix a
Lie algebra $\mathfrak{g}$.  We want to use the
Campbell-Baker-Hausdorff formula, where the formal symbols $a$ and $b$
in the series $C(a,b)$ are now replaced by elements of $\mathfrak{g}$,
and we want to consider the resulting Lie series as an element of
$\mathfrak{g}$ rather than as a formal series.  In other words, we
want to evaluate the brackets within the Lie algebra $\mathfrak{g}$,
but in general, the series might not be well defined in $\mathfrak{g}$
since it will typically contain infinitely many nonzero terms.
However, let us introduce commuting formal variables $s$ and $t$
commuting with $\mathfrak{g}$.  We consider the
Campbell-Baker-Hausdorff series $C(sg_1,tg_2)$ for any $g_1, g_2 \in
\mathfrak{g}$.  {}From (\ref{real CH formula}), (\ref{C decomp 1}) and
(\ref{C decomp 2}), we see that each $H_{j,k}$ for $j,k \in
\mathbb{N}$ in the Lie series $C(sg_1,tg_2)$ involves only finitely
many brackets in $\mathfrak{g}[[s,t]]$, and therefore when we evaluate
brackets in $\mathfrak{g}$, the coefficient of a given $s^j t^k$ in
$C(sg_1,tg_2)$ is well defined in $\mathfrak{g}$.  Thus $C(sg_1,tg_2)$
is a well defined element of $\mathfrak{g}[[s,t]]$.  Furthermore, note
that if, more generally, $g_1, g_2 \in \mathfrak{g}[[s,t]]$, then for
$j,k \in \mathbb{N}$, each $H_{j,k}$ in the Lie series $C(sg_1,tg_2)$
is a sum of terms
of degree greater than or equal to $j$ in $s$ and $k$ in $t$, so
that $C(sg_1,tg_2)$ is still well defined in $\mathfrak{g}[[s,t]]$.
Of course the exponentials $e^{sg_1}$, $e^{tg_2}$ and $e^{C(sg_1,
tg_2)}$ are elements of the algebra $U(\mathfrak{g})[[s,t]]$ of formal
power series over the universal enveloping algebra $U(\mathfrak{g})$.

Fix a vector space decomposition of the Lie algebra $\mathfrak{g}$:
\begin{equation}\label{decomposition}
\mathfrak{g} = \mathfrak{g}^- \oplus \mathfrak{g}^+ ,
\end{equation}  
so that $\mathfrak{g} [[s,t]] = \mathfrak{g}^- [[s,t]] \oplus
\mathfrak{g}^+ [[s,t]]$.  Let
\[\pi^\pm : \mathfrak{g} [[s,t]] \longrightarrow \mathfrak{g}^\pm
[[s,t]] \] 
be the corresponding projection maps.  Then we have:

\begin{thm}[The Factorization Theorem]\label{main theorem}
The map
\begin{equation}\label{mapC}
C: s\mathfrak{g}^-[[s,t]] \times t\mathfrak{g}^+[[s,t]]
\longrightarrow s\mathfrak{g}^-[[s,t]] \oplus t\mathfrak{g}^+[[s,t]] 
\end{equation}
defined by the Campbell-Baker-Hausdorff formula and given by 
\begin{equation}\label{main formula}
e^{sg^-} e^{tg^+} = e^{C(sg^-,tg^+)}
\end{equation}
for $g^\pm \in \mathfrak{g}^\pm[[s,t]]$ is a bijection.  The lowest 
order terms of $C(sg^-,tg^+)$ are described as follows:
\begin{equation}\label{CBH conditions}
C(sg^-,tg^+) = sg^- + tg^+ + st \frac{1}{2} [g^-, g^+] + st q(s,t),
\end{equation}
where $q(s,t) \in \mathfrak{g}[[s,t]]$ and $q (0,0) = 0$.  Moreover, 
for $h^\pm \in \mathfrak{g}^\pm[[s,t]]$, the lowest order terms of
\[C^{-1}(sh^- + th^+) = (sg^-, tg^+)\] 
are described as follows:
\begin{eqnarray}
g^- = h^- + \frac{1}{2} t \pi^- ([h^+, h^-]) + tr^-(s,t),
\label{conditions 1} \\ 
g^+ = h^+ + \frac{1}{2} s \pi^+ ([h^+, h^-]) + sr^+(s,t), 
\label{conditions 2}
\end{eqnarray}
where $r^\pm (s,t) \in \mathfrak{g}^\pm [[s,t]]$ and $r^\pm (0,0) =
0$.
\end{thm}

{\it Proof:} \hspace{.2cm} For $g^\pm \in \mathfrak{g}^\pm [[s,t]]$,
the expression $e^{sg^-} e^{tg^+}$ is well defined in
$U(\mathfrak{g})[[s,t]]$, and by the Campbell-Baker-Hausdorff theorem
and the discussion above, there exists a unique element $C(sg^-,tg^+) 
\in \mathfrak{g}[[s,t]]$ such that (\ref{main formula}) holds; 
moreover, (\ref{CBH conditions}) also holds.

For $g^\pm \in \mathfrak{g}^\pm [[s,t]]$, write $g^\pm = \sum_{m \in
\mathbb{N}} g^\pm_m$ where $g^\pm_m$ is homogeneous of degree $m$ in
$s$, and write
\begin{equation}\label{g pm sub m}
g^\pm_m = \sum_{n \in \mathbb{N}} g^\pm_{m,n}
\end{equation}
where $g^\pm_{m,n} \in \mathfrak{g}^\pm [s,t]$ is homogeneous of
degree $m$ in $s$ and degree $n$ in $t$, i.e., 
$g^\pm_{m,n}s^{-m}t^{-n} \in \mathfrak{g}^\pm$.  To prove that $C$ is
bijective, given $h^\pm \in \mathfrak{g}^\pm [[s,t]]$, we will use
recursion on $m$ and $n$ to construct unique $g^\pm \in
\mathfrak{g}^\pm [[s,t]]$ such that 
\begin{equation}\label{maineq}
C(sg^-, tg^+) = sh^- + th^+.
\end{equation}
We will use the notation $h^\pm_{m}$ and $h^\pm_{m,n}$, by analogy
with (\ref{g pm sub m}).

For $g^\pm \in \mathfrak{g}^\pm[[s,t]]$, consider the series in
$\mathfrak{g}[[s,t]]$ given by the Campbell-Baker-Hausdorff formula
\begin{equation}\label{C=H}
C(sg^-, tg^+) = \sum_{j,k \in \mathbb{N}} H_{j,k} (sg^-, tg^+)
\end{equation}
where $H_{j,k} (sg^-, tg^+)$ involves brackets containing exactly
$j$ of the elements $sg^-$ and exactly $k$ 
of the elements $tg^+$, as in (\ref{C decomp 1}) and (\ref{C decomp 2}).
In particular, $H_{j,k} (sg^-, tg^+)$ is 
a sum of terms of degree greater than or
equal to $j$ in $s$ and degree greater than or equal to $k$ in $t$.
{}From (\ref{real CH formula}),
\begin{equation}\label{degrees of s}
C(sg^-, tg^+) = tg^+ + \left( \frac{\mathrm{ad} \;
tg^+}{e^{\mathrm{ad} \; tg^+} - 1} \right) (sg^-) + p(s,t)
\end{equation}
where $p(s,t)$ is a Lie series in $sg^-$ and $tg^+$ whose terms have
degree greater than one in $sg^-$ and degree greater than zero in
$tg^+$, and in particular, whose terms have 
degree greater than one in $s$ and degree
greater than zero in $t$.  We set $C(sg^-, tg^+) = sh^- + th^+$ in
(\ref{degrees of s}).  Equating the terms of degree zero in $s$ is
equivalent to the ``initial conditions''
\begin{equation}\label{initial plus}
g^+_0 = \sum_{n \in \mathbb{N}} g^+_{0,n} = \sum_{n \in \mathbb{N}} 
h^+_{0,n} = h^+_0 .
\end{equation} 
Moreover, equating the terms of degree zero in $t$ is equivalent to 
the initial conditions
\begin{equation}\label{initial minus}
\sum_{m \in \mathbb{N}} g^-_{m,0} = \sum_{m \in \mathbb{N}} 
h^-_{m,0}.
\end{equation}
 
Equating the terms of degree one in $s$ in equation (\ref{degrees of
s}) and using the initial conditions (\ref{initial plus}) 
amounts to the formula
\begin{eqnarray*} 
sh^-_0 + th^+_1 &=& t g^+_1 + \left( \frac{\mathrm{ad} \;
th^+_0}{e^{\mathrm{ad}\; th^+_0} - 1} \right) \left( s g^-_0 \right)
\\ 
&=& t g^+_1 + \left( \sum_{k \in \mathbb{N}} \frac{B_k}{k !}
(\mathrm{ad} \; th^+_0)^k \right) \left( s g^-_0 \right),
\end{eqnarray*}
where $B_k$, $k \in \mathbb{N}$, are the Bernoulli numbers, defined by
the generating function
\[\sum_{k \in \mathbb{N}} \frac{B_k}{k!} x^k = \frac{x}{e^x - 1} . \]
Equivalently (using the decomposition (\ref{decomposition})),
\begin{eqnarray}
sg^-_0 &=& sh^-_0 - s\pi^- \left( \sum_{k > 0}
\frac{B_k}{k!}  (\mathrm{ad} \; th^+_0)^k g^-_0 \right),
\label{one in s minus}\\ 
tg^+_1 &=& th^+_1 - s\pi^+ \left( \sum_{k > 0} \frac{B_k}{k!}
(\mathrm{ad} \; th^+_0)^k g^-_0 \right) \label{one in s plus}.
\end{eqnarray}
Equating the terms of degree one in $s$ and one in $t$ in
equations (\ref{one in s minus}) and (\ref{one in s plus}) and 
using the initial conditions
(\ref{initial minus}) is equivalent to the following information:
\begin{eqnarray}\label{g01-}
g^-_{0,1} &=& h^-_{0,1} - B_1 \pi^- ( [th^+_{0,0}, h^-_{0,0} ]) 
\nonumber \\
&=& h^-_{0,1} + \frac{1}{2} t \pi^- ( [h^+_{0,0}, h^-_{0,0} ]) 
\end{eqnarray}
and
\begin{eqnarray}\label{g10+}
g^+_{1,0} &=& h^+_{1,0} - s B_1 \pi^+ ( [h^+_{0,0}, h^-_{0,0} ])
\nonumber \\
&=& h^+_{1,0} + \frac{1}{2} s \pi^+ ([h^+_{0,0}, h^-_{0,0} ]) .
\end{eqnarray}
The conditions (\ref{initial plus}), (\ref{initial minus}),
(\ref{g01-}) and (\ref{g10+}) are together equivalent to
(\ref{conditions 1}) and (\ref{conditions 2}).

We will use recursion on the subscripts to construct and uniquely
determine all the $g^\pm_{m,n}$ by equating the coefficients of
appropriate powers of $s$ and $t$ in (\ref{maineq}).  So far, we have
the following: Equating the coefficients of $s^{0}t^{n}$ ($n \ge 1$)
in (\ref{maineq}) is equivalent to the information $g^+_{0,n-1} =
h^+_{0,n-1}$ (\ref{initial plus}); the equation for $s^{m}t^{0}$ ($m
\ge 1$) is equivalent to $g^-_{m-1,0} = h^-_{m-1,0}$ (\ref{initial
minus}); and the equation for $s^{1}t^{1}$ is equivalent to
(\ref{g01-}) and (\ref{g10+}), using the special cases $g^+_{0,0} =
h^+_{0,0}$ and $g^-_{0,0} = h^-_{0,0}$ of (\ref{initial plus}) and
(\ref{initial minus}).  Note that (\ref{one in s minus}) and (\ref{one
in s plus}) do not serve to uniquely construct $g^-_{0,n}$ for $n > 1$
or $g^+_{m,0}$ for $m > 1$ as desired, since we must still use the
recursive procedure below to uniquely express the components of
$g^-_{0}$ on the right-hand sides of (\ref{one in s minus}) and
(\ref{one in s plus}) in terms of $h$'s.  (The general recursion below
will redo the cases $g^-_{0,1}$ and $g^+_{1,0}$, but we needed the
precise formulas (\ref{g01-}) and (\ref{g10+}).)

Let $m,n > 0$.  Equating the coefficients of $s^{m}t^{n}$ in
(\ref{maineq}) is equivalent to the equation 
\begin{equation}\label{generalstep}
sg^-_{m-1,n} + tg^+_{m,n-1} + l_{m,n} = sh^-_{m-1,n}+ th^+_{m,n-1},
\end{equation}
where $l_{m,n}$ is an explicit linear combination, homogeneous of
degree $m$ in $s$ and of degree $n$ in $t$, of brackets of two or more
elements of the form $sg^-_{p,q}$ and $tg^+_{p,q}$ (with at least one
of each of these two types) with $p < m$ and $q < n$.  We
equivalently have
\[ sg^-_{m-1,n} = sh^-_{m-1,n} - \pi^- (l_{m,n}) \]
and
\[ tg^+_{m,n-1} = th^+_{m,n-1} - \pi^+ (l_{m,n}). \]
Proceeding through the sequence (for example)
\[ (m,n) = (1,1);(1,2),(2,1);(1,3),(2,2),(3,1);(1,4),\dots , \]
we see that we have an effective recursive procedure for uniquely
computing all the $g^\pm_{m,n}$ in terms of the $h^\pm_{m,n}$, Lie
brackets, and the projections $\pi^{\pm}$.  In particular, the map $C$
is bijective. \pfbox

\begin{prob}\label{problem}
{\em The map $C^{-1}$ is given by a precise recursive procedure.  We
propose the following problem:  Find a closed form for this map.
See also Problem \ref{problem2} below.}
\end{prob}     

\begin{rema}\label{several subspaces}
{\em Theorem \ref{main theorem} and our method of proof, based on the
decomposition (\ref{decomposition}), generalize to the case of a
decomposition of $\mathfrak{g}$ into an arbitrary finite number of
subspaces.}
\end{rema}     

The proof of Theorem \ref{main theorem} yields the following more
precise information about how the elements $g^\pm$ are built {}from
the elements $h^\pm$ using commutators and the projections
$\pi^{\pm}$, under the assumption (which we will remove in Theorem
\ref{generalwordsthm} below) that $h^\pm \in \mathfrak{g}$, i.e., that
the elements $h^\pm$ do not involve $s$ or $t$:

\begin{prop}\label{wordsprop}
In the setting of Theorem \ref{main theorem}, suppose that $h^\pm \in
\mathfrak{g}$.  Write $\mathbb{F}\langle h^\pm ; \pi^\pm \rangle$
(resp., $\mathbb{F}\langle sh^-,th^+ ; \pi^\pm \rangle$) for the
smallest Lie subalgebra of $\mathfrak{g}$ (resp.,
$\mathfrak{g}[[s,t]]$) containing the elements $h^-$ and $h^+$ (resp.,
$sh^-$ and $th^+$) and closed under the projections $\pi^\pm$, i.e.,
compatible with the decomposition (\ref{decomposition}) in the sense
that it is the direct sum of its intersections with $\mathfrak{g}^\pm$
(resp., $\mathfrak{g}^\pm [[s,t]]$).  We have:

(a) The Lie algebra $\mathbb{F}\langle sh^-,th^+ ; \pi^\pm \rangle$ is
$\mathbb{N} \times \mathbb{N}$-graded by means of the decomposition
\begin{equation}
\mathbb{F}\langle sh^-,th^+ ; \pi^\pm \rangle =
\coprod_{m,n \in \mathbb{N}}\mathbb{F}\langle h^\pm ; \pi^\pm
\rangle_{m,n}s^m t^n,
\end{equation}
where $\mathbb{F}\langle h^\pm ; \pi^\pm\rangle_{m,n}$ is the subspace
of $\mathbb{F}\langle h^\pm ; \pi^\pm\rangle$ spanned by the elements
built {}from commutators involving exactly $m$ elements $h^-$ and
exactly $n$ elements $h^+$, and {}from the projections $\pi^\pm$.
(Warning: The subspaces $\mathbb{F}\langle h^\pm ;
\pi^\pm\rangle_{m,n}$ might not be disjoint; for instance, we might
have, say, $[h^+,h^-]=h^+$.)

(b) Consider the formal completion
\begin{equation}\label{formalcompletion}
\mathbb{F}\langle \langle sh^-,th^+ ; \pi^\pm \rangle \rangle =
\prod_{m,n \in \mathbb{N}}\mathbb{F}\langle h^\pm ; \pi^\pm
\rangle_{m,n}s^m t^n
\end{equation}
of $\mathbb{F}\langle sh^-,th^+ ; \pi^\pm \rangle$ in
$\mathfrak{g} [[s,t]]$, so that $\mathbb{F}\langle \langle
sh^-,th^+ ; \pi^\pm \rangle \rangle$ is naturally a Lie subalgebra of
$\mathbb{F}\langle h^\pm ; \pi^\pm \rangle[[s,t]]$ stable under
$\pi^\pm$.  We have:
\begin{equation}
sg^-,tg^+ \in \mathbb{F}\langle \langle sh^-,th^+ ; \pi^\pm \rangle \rangle.
\end{equation}
\end{prop}

{\it Proof:} \hspace{.2cm} Part (a) is clear.  To prove (b), it is
sufficient to show that
\begin{equation}\label{formulafor(b)}
(sg^-)_{m,n},(tg^+)_{m,n} \in \mathbb{F}\langle h^\pm ; \pi^\pm
\rangle_{m,n}s^m t^n
\end{equation}
for all $m$ and $n$ (using the notation (\ref{g pm sub m}) for
homogeneous components of elements of $\mathfrak{g}[[s,t]]$).  We
proceed through the proof of Theorem \ref{main theorem} and indicate
the special information that we have in this situation.  Formula
(\ref{maineq}) remains the same, but since the elements $h^{\pm}$ do
not involve $s$ or $t$, we do not need to consider the components
$h^\pm_{m}$ or $h^\pm_{m,n}$.  {}From (\ref{maineq}) we find that
(\ref{initial plus}) and (\ref{initial minus}) become, respectively:
\begin{eqnarray}
tg^+_0 &=& th^+, \label{m=0} \\
s\sum_{m \in \mathbb{N}} g^-_{m,0} &=& sh^-, \label{n=0}
\end{eqnarray}
so that $g^+_0=h^+$ is independent of $s$ and $t$, and $\sum
g^-_{m,0}=g^-_{0,0}=h^-$ and is also independent of $s$ and $t$.
Also, (\ref{g01-}) and (\ref{g10+}) become, respectively:
\begin{eqnarray}
sg^-_{0,1} &=& \frac{1}{2} \pi^- ([th^+,sh^-]), \label{sg01} \\
tg^+_{1,0} &=& \frac{1}{2} \pi^+ ([th^+,sh^-]). \label{tg10}
\end{eqnarray}
For $m,n > 0$, the right-hand side of (\ref{generalstep}) is 0, and
(\ref{generalstep}) becomes:
\begin{equation}
sg^-_{m-1,n} + tg^+_{m,n-1} + l_{m,n} = 0,
\end{equation}
where $l_{m,n}$ is an (explicit) linear combination as indicated in
the proof.  Now we just use the inductive procedure described in the
proof to establish (\ref{formulafor(b)}) by induction on $(m,n)$.  The
cases $m=0$ and $n=0$ are covered by (\ref{m=0}) and (\ref{n=0}),
respectively, and the case $(1,1)$ follows {}from (\ref{sg01}) and
(\ref{tg10}); $l_{1,1}$ is a multiple of $s^1 t^1 [h^+,h^-]$.  The
general inductive step is clear, and the result is proved.  \pfbox

\begin{rema}\label{canonicalfactorizationseries}
{\em The proof constructs $sg^-$ and $tg^+$, using the
Campbell-Baker-Hausdorff series, as canonical formal series of
``words'' involving brackets of $sh^-$ and $th^+$, and the projections
$\pi^\pm$, independently of the Lie algebra $\mathfrak{g}$ or of
$\pi^\pm$ or of the elements $h^\pm$.  Let us call these two formal
series of words the {\it canonical factorization series} and let us
write them as
\begin{equation}
F^\pm (sh^-,th^+;\pi^\pm).
\end{equation}
They are analogues, in some sense, of the Campbell-Baker-Hausdorff
series.}
\end{rema}

Now we remove the assumption $h^\pm \in \mathfrak{g}$ in Proposition
\ref{wordsprop}, using Propsition \ref{wordsprop} to obtain the
corresponding information in the general case.  We write $\mathbb{F}
\langle sh^-,th^+ ; \pi^\pm \rangle$ for the smallest Lie subalgebra
of $\mathfrak{g}[[s,t]]$ containing $sh^-$ and $th^+$ and closed under
$\pi^\pm$ and we define the formal completion
$\mathbb{F}\langle\langle sh^-,th^+ ; \pi^\pm \rangle\rangle$ of
$\mathbb{F}\langle sh^-,th^+ ; \pi^\pm \rangle$ to be the vector space
of formal (possibly infinite) linear combinations of ``words''
involving brackets of the elements $sh^-$ and $th^+$, and the
projections $\pi^\pm$.  This space is well defined because there are
only finitely many words involving $s^m t^n$ for fixed $m,n \in
\mathbb{N}$.  The space $\mathbb{F}\langle\langle sh^-,th^+ ; 
\pi^\pm \rangle\rangle$ is clearly a Lie subalgebra of 
$\mathfrak{g}[[s,t]]$ stable under $\pi^\pm$.  Note that this Lie
algebra is an analogue of the Lie algebra $\mathbb{F}\langle\langle
a,b \rangle\rangle$ in Section 2.  In the special case that $h^\pm \in
\mathfrak{g}$, this Lie algebra agrees with the already-defined Lie
algebra (\ref{formalcompletion}).  We have the following
generalization of Proposition \ref{wordsprop}:

\begin{thm}\label{generalwordsthm}
In the general setting of Theorem \ref{main theorem} (in the absence
of the assumption $h^\pm \in \mathfrak{g}$), we have:
\begin{equation}
sg^-,tg^+ \in \mathbb{F}\langle \langle sh^-,th^+ ; \pi^\pm \rangle \rangle,
\end{equation}
and $sg^-$ and $tg^+$ are given by the canonical factorization series
(recall Remark \ref{canonicalfactorizationseries}):
\begin{equation}
sg^- = F^- (sh^-,th^+;\pi^\pm), \;\;\; tg^+ = F^+ (sh^-,th^+;\pi^\pm).
\end{equation}
Moreover, the lowest-order terms in $sg^-$ and $tg^+$ are given by:
\begin{eqnarray}
sg^- = sh^- + \frac{1}{2}\pi^-([th^+,sh^-]) + u^-(s,t), \label{lowest order F -}\\
tg^+ = th^+ + \frac{1}{2}\pi^+([th^+,sh^-]) + u^+(s,t), \label{lowest order F +}
\end{eqnarray}
where $u^\pm(s,t) \in \mathbb{F}\langle \langle sh^-,th^+ ; \pi^\pm
\rangle \rangle^{\pm}$ are formal (possibly infinite) linear
combinations of words involving at least three occurrences of 
$sh^-$ and $th^+$ (including at least one of each).
\end{thm}

{\it Proof:} \hspace{.2cm} Write $\mathfrak{h}$ for the Lie algebra
$\mathfrak{g}[[s,t]]$ and apply Proposition \ref{wordsprop} to the Lie
algebra $\mathfrak{h}$ in place of $\mathfrak{g}$ and
$\mathfrak{h}[[s_1,t_1]]$ in place of $\mathfrak{g}[[s,t]]$, with
$s_1$ and $t_1$ new formal variables.  We find that given our elements
$h^\pm \in \mathfrak{g}^{\pm}[[s,t]]$, we have that the formula
\begin{equation}\label{mainformula_1}
e^{s_1 g^{-}_{1}}e^{t_1 g^{+}_{1}} = e^{s_1 h^{-} + t_1 h^{+}}
\end{equation}
determines unique elements (by Theorem \ref{main theorem})
\begin{equation}
g_{1}^{\pm} \in \mathfrak{h}[[s_1,t_1]],
\end{equation}
and by Proposition \ref{wordsprop}, for all $m,n \in \mathbb{N}$, the
coefficient of $s_{1}^m t_{1}^n$ in $s_1 g^{-}_{1}$ and in $t_1
g^{+}_{1}$ lies in $\mathbb{F} \langle h^\pm ; \pi^\pm \rangle_{m,n}$,
where $\mathbb{F} \langle h^\pm ; \pi^\pm \rangle$ and $\mathbb{F}
\langle h^\pm ; \pi^\pm \rangle_{m,n}$ are defined as in Proposition
\ref{wordsprop} in the present case.  Moreover, $s_1 g^{-}_{1}$ and
$t_1 g^{+}_{1}$ are expressed by the canonical factorization series in
terms of $s_1 h^{-}$ and $t_1 h^{+}$, and $\pi^{\pm}$, and the
low-order terms with respect to $s_1$ and $t_1$ in $g_{1}^{\pm}$ are
given by (\ref{conditions 1}) and (\ref{conditions 2}) with $s_1$ and
$t_1$ in place of $s$ and $t$ and with $\mathbb{F} \langle h^\pm ;
\pi^\pm \rangle$ in place of $\mathfrak{g}$.  We may set $s=s_1$ and
$t=t_1$ in (\ref{mainformula_1}), and we see that the elements
$g^{\pm}$ are determined {}from the elements
\begin{equation}
g_{1}^{\pm} \in \mathbb{F} \langle h^\pm ; \pi^\pm \rangle [[s_1,t_1]]
\end{equation}
(which are uniquely determined by the formula (\ref{mainformula_1})),
by the specialization
\begin{equation}
g^{\pm} = g_{1}^{\pm}\vert_{s_{1}=s,\;t_{1}=t}.
\end{equation}
Moreover,
\begin{eqnarray}
s_{1}g_{1}^- = s_{1}h^- + \frac{1}{2} s_{1}t_{1}\pi^-([h^+,h^-]) + 
s_{1}t_{1}r_{1}^-(s_{1},t_{1}), \\
t_{1}g_{1}^+ = t_{1}h^+ + \frac{1}{2} s_{1}t_{1}\pi^+([h^+,h^-]) + 
s_{1}t_{1}r_{1}^+(s_{1},t_{1}),
\end{eqnarray}
where
\begin{equation}
r_{1}^{\pm}(s_{1},t_{1}) \in \mathbb{F} \langle h^\pm ; \pi^\pm
\rangle^{\pm}[[s_1,t_1]], \;\;\; r_{1}^{\pm}(0,0) = 0,
\end{equation}
and for all $m,n \ge 1$ with either $m$ or $n \ge 2$, the coefficient
of $s_{1}^m t_{1}^n$ in $s_1 t_1 r_{1}^{\pm}(s_{1},t_{1})$ is a linear
combination of elements of $\mathbb{F} \langle h^\pm ; \pi^\pm
\rangle$ built according to the canonical factorization series {}from
commutators involving exactly $m$ elements $h^-$ and exactly $n$
elements $h^+$, and {}from the projections $\pi^{\pm}$.  This gives the
desired result.  \pfbox

\begin{prob}\label{problem2}
{\em Find closed forms for the canonical factorization series
$F^{\pm}$.  (Cf. Problem \ref{problem} and the classical formula
(\ref{real CH formula}) for the Campbell-Baker-Hausdorff series.)}
\end{prob}

\begin{rema}\label{directuniqueness}
{\em Here we give an alternate, direct, simple proof of the uniqueness
of the factors in the product $e^{sg^-} e^{tg^+}$ in (\ref{main
formula}) (i.e., the injectivity of the map $C$ in Theorem \ref{main
theorem}), under the extra hypothesis that the subspaces
$\mathfrak{g}^\pm$ are Lie subalgebras (but see the next remark for
the removal of this extra hypothesis).  The following argument also
works more generally for the analogous uniqueness when the Lie algebra
$\mathfrak{g}$ is given as a finite direct sum of any number of Lie
subalgebras.  We use the Poincar\'e-Birkhoff-Witt theorem.  Write
$$P:U(\mathfrak{g})\longrightarrow U(\mathfrak{g}^-)$$ for the
projection with respect to the decomposition $$U(\mathfrak{g}) =
U(\mathfrak{g}^-) \oplus U(\mathfrak{g}^-)
U(\mathfrak{g}^+)\mathfrak{g}^+ = U(\mathfrak{g}^-)
\oplus U(\mathfrak{g})\mathfrak{g}^+,$$
coming {}from the decomposition
$$U(\mathfrak{g}) = U(\mathfrak{g}^-) \otimes U(\mathfrak{g}^+),$$
and extend $P$ canonically to $U(\mathfrak{g})[[s,t]]$.  Given
\begin{equation}\label{ee=ee}
e^{sg^-} e^{tg^+} = e^{sg^{-}_{1}} e^{tg^{+}_{1}}
\end{equation}
($g^{\pm}_{1} \in \mathfrak{g}^{\pm}[[s,t]]$), simply apply $P$ to get
$e^{sg^-} = e^{sg^{-}_{1}}$ and hence $g^- = g^{-}_{1}$, and {}from
this, $g^+ = g^{+}_{1}$.}
\end{rema}     

\begin{rema}\label{generaldirectuniqueness}
{\em Here we remove the extra hypothesis in Remark
\ref{directuniqueness} that the subspaces $\mathfrak{g}^\pm$ be Lie
subalgebras, proving the uniqueness in Theorem \ref{main theorem} in
general; as in Remark \ref{directuniqueness}, this argument works in
the general setting of Remark \ref{several subspaces}.  Let 
$$\lambda : S(\mathfrak{g}) \longrightarrow U(\mathfrak{g})$$ 
be the standard symmetrization map---a linear isomorphism by the
Poincar\'e-Birkhoff-Witt theorem---determined by:
$$\lambda (g_1 \cdots g_k) = \frac{1}{k!}\sum_{\sigma} g_{\sigma(1)} \cdots
g_{\sigma(k)},$$
where $k \ge 0$, $g_1,\dots,g_k \in \mathfrak{g}$ and $\sigma$ ranges
through the symmetric group on $k$ letters.  Then by the
Poincar\'e-Birkhoff-Witt theorem,
\begin{eqnarray*}
U(\mathfrak{g}) &=& \lambda(S(\mathfrak{g}^-)) \otimes 
\lambda(S(\mathfrak{g}^+))
\nonumber \\
&=& \lambda(S(\mathfrak{g}^-)) \oplus \lambda(S(\mathfrak{g}^-))
\lambda(\mathfrak{g}^{+}S(\mathfrak{g}^+)).
\end{eqnarray*}
Let
$$P:U(\mathfrak{g})\longrightarrow \lambda(S(\mathfrak{g}^-))$$
be the corresponding projection, and extend $P$ canonically to the
projection
$$P:U(\mathfrak{g})[[s,t]]\longrightarrow \lambda(S(\mathfrak{g}^-))
[[s,t]].$$
Now with $g^- \in \mathfrak{g}^{-}[[s,t]]$ as above, we have
$$e^{sg^-} \in \lambda(S(\mathfrak{g}^-))[[s,t]]$$
(and similarly for $tg^+$).  Indeed, the coefficient of each monomial
in $s$ and $t$ in $e^{sg^-}$ coincides with the coefficient of the
same  monomial in a suitable {\it finite} linear combination of powers
of $sg^-$, and for any $k \ge 0$,
$$(sg^-)^k \in \lambda(S(\mathfrak{g}^-))[[s,t]],$$
since the map $\lambda$ extends canonically to the natural map
$$\lambda : S(\mathfrak{g})[[s,t]] \longrightarrow U(\mathfrak{g})[[s,t]],$$ 
and $(sg^-)^k$ is the image of $(sg^-)^k$ viewed as an element of 
$S(\mathfrak{g}^-)[[s,t]]$.  Given (\ref{ee=ee}), we can now apply $P$
just as in Remark \ref{directuniqueness} to get $e^{sg^-} =
e^{sg^{-}_{1}}$, giving the uniqueness.}
\end{rema}     

\begin{rema}\label{universal vs lie}
{\em Recall that the nontrivial part of the Campbell-Baker-Hausdorff
theorem is that the element $C(a,b)$ (see formula (\ref{CH})) of the
universal enveloping algebra of the free Lie algebra over $a$ and $b$
is in fact a Lie series.  In the proof of bijectivity in Theorem
\ref{main theorem}, in general, the operations take place in
$U(\mathfrak{g})[[s,t]]$; however, the projections $\pi^\pm$ are not
defined on $U(\mathfrak{g})[[s,t]]$.  Rather, it is necessary to use
the fact that $C(sg^-,tg^+)$ is a Lie series, to take brackets in the
Lie algebra $\mathfrak{g} [[s,t]]$, and then to project using
$\pi^\pm$.  One cannot use (associative) words in $g^\pm$ to obtain
$C^{-1}(sh^- + th^+)$, in contrast with the situation for
$C(sg^-,tg^+)$.  However, Theorem \ref{generalwordsthm} shows that
$C^{-1}(sh^- + th^+)$ can be obtained using words in $sh^-$ and $th^+$ and
the canonical projections $\pi^\pm$, i.e., the correct setting for
$C^{-1}(sh^- + th^+)$ is $\mathbb{F}\langle \langle sh^-, th^+; 
\pi^\pm \rangle \rangle$, as opposed to $U(\mathfrak{g})[[s,t]]$ for
$C(sg^-,tg^+)$.}
\end{rema}     

\begin{rema}\label{method of proof}
{\em The method used to prove that the map $C$ in Theorem \ref{main
theorem} is bijective is similar to the method used to prove the
``formal uniformization'' result Theorem 2.2.4 in \cite{H3} for a
certain representation of the Virasoro algebra, given by $L_n \mapsto
- x^{n+1} \frac{d}{dx} \in \mathrm{End} (\mathbb{C}[x,x^{-1}])$ and $c
= 0$, and to prove the analogous result, Theorem 2.3.4, in \cite{B1},
for an analogous representation of the $N = 1$ Neveu-Schwarz algebra.
The method is similar in that once one has an appropriate equation
involving formal series and one has appropriate projections, one can
solve the equation recursively.  However, the settings in the proofs
for these two cases are very different {}from the setting in the present
proof: In \cite{H3} and \cite{B1} the proofs take place in certain
formal function algebras, whereas in the present paper the proof takes
place in a certain Lie algebra. Note the subtle issue that in a
universal enveloping algebra or formal extension there are no
appropriate projections available for the use of this method (recall
Remark \ref{universal vs lie}).  In fact a crucial observation in this
paper is that even though there is no such setting involving universal
enveloping algebras, one can in fact find a very general setting to
which the method can be applied, a setting very different {}from the
ones in \cite{H3} and \cite{B1}.  As a benefit, we are able to obtain
the factorization results in \cite{H3} and \cite{B1} in a uniform,
simple way, without the need to pass to a central extension (see
Applications \ref{virasoro example} and \ref{N = 1} below, Theorem
2.2.4, Proposition 4.2.1 and Corollary 4.2.2 in \cite{H3}, and Theorem
2.3.4, Proposition 2.6.1 and Corollary 2.6.2 in \cite{B1}).}
\end{rema}

\begin{cor}[Formal Algebraic Uniformization]\label{switching corollary}
There exists a unique bijection 
\begin{equation}
\Psi: t\mathfrak{g}^+[[s,t]] \times s\mathfrak{g}^-[[s,t]]
\longrightarrow s\mathfrak{g}^-[[s,t]] \times t\mathfrak{g}^+[[s,t]] 
\end{equation}
such that for $g^\pm \in \mathfrak{g}^\pm[[s,t]]$,
\begin{equation}\label{switching equation}
e^{tg^+} e^{sg^-} = e^{s\Psi^-} e^{ t\Psi^+}
\end{equation}
with $s\Psi^- = \pi^- \circ \Psi (tg^+, sg^-)$ and $t \Psi^+ = \pi^+
\circ \Psi (tg^+, sg^-)$.  Moreover, the lowest order terms of the
$\Psi^\pm$ are described as follows:
\begin{eqnarray}
\Psi^- = g^- + t \pi^-([g^+, g^-]) + t{\cal P}^-(s,t) \label{psi
conditions 1}\\ 
\Psi^+ = g^+ + s \pi^+([g^+, g^-]) + s{\cal P}^+(s,t) \label{psi
conditions 2} 
\end{eqnarray}
where ${\cal P}^\pm (s,t) \in \mathfrak{g}^\pm [[s,t]]$ and ${\cal
P}^\pm (0,0) = 0$.
\end{cor}

{\it Proof:} \hspace{.2cm} Let $C'$ be the analogue of the map $C$
(see (\ref{mapC}), (\ref{main formula})) with the roles of $s$ and $t$
and of $\mathfrak{g}^-$ and $\mathfrak{g}^+$ reversed, and let
$\sigma$ be the isomorphism given by
\begin{eqnarray*} 
\sigma: t\mathfrak{g}^+[[s,t]] \oplus s\mathfrak{g}^-[[s,t]]
&\longrightarrow& s\mathfrak{g}^-[[s,t]] \oplus t\mathfrak{g}^+[[s,t]] \\
(th^+, sh^-) &\mapsto& (sh^-, th^+) . 
\end{eqnarray*} 
Then 
\[\Psi = C^{-1} \circ \sigma \circ C' : t\mathfrak{g}^+[[s,t]] \times 
s\mathfrak{g}^-[[s,t]] \longrightarrow s\mathfrak{g}^-[[s,t]] \times
t\mathfrak{g}^+[[s,t]]\] 
is a bijection satisfying (\ref{switching equation}).  The conditions
given by equations (\ref{psi conditions 1}) and (\ref{psi conditions
2}) for the lowest order terms follow {}from equations (\ref{CBH
conditions}), (\ref{conditions 1}), and (\ref{conditions 2}). \pfbox

Theorem \ref{generalwordsthm} along with the Campbell-Baker-Hausdorff
theorem of course gives the corresponding additional information about
the elements $\Psi^{\pm}$ in Corollary \ref{switching corollary}:

\begin{cor}\label{switchingwordscorollary}
We have:
\begin{equation}
s\Psi^-,t\Psi^+ \in \mathbb{F}\langle \langle sg^-,tg^+ ;
\pi^\pm \rangle \rangle
\end{equation}
and $s\Psi^-$ and $t\Psi^+$ are given by canonical formal series,
which we write as:
\begin{equation}
s\Psi^- = G^- (sg^-,tg^+;\pi^\pm), \;\;\; t\Psi^+ = G^+ (sg^-,tg^+;\pi^\pm).
\end{equation}
The lowest-order terms are given by:
\begin{eqnarray}
G^- (sg^-,tg^+;\pi^\pm) = sg^- + \pi^-([tg^+,sg^-]) + v^-(s,t), 
\label{lowest order G -}\\
G^+ (sg^-,tg^+;\pi^\pm) = tg^+ + \pi^+([tg^+,sg^-]) + v^+(s,t),  
\label{lowest order G +}
\end{eqnarray}
where $v^\pm(s,t) \in \mathbb{F}\langle \langle sg^-,tg^+ ; \pi^\pm
\rangle \rangle^{\pm}$ are formal (possibly infinite) linear
combinations of words involving at least three occurrences of 
$sg^-$ and $tg^+$ (including at least one of each).  \pfbox
\end{cor}

\begin{prob}\label{problem3}
{\em Let us call the two series $G^{\pm}$ the {\it formal algebraic
uniformization series}.  They are essentially compositions, in a
suitable sense, of the Campbell-Baker-Hausdorff series and the
canonical factorization series, incorporating the twist $\sigma$ in
the proof of Corollary \ref{switching corollary}.  Find closed forms
for the formal algebraic uniformization series.  (Cf. Problem
\ref{problem2}.)}
\end{prob}

\section{Applications to $\mathbb{Z}$-graded Lie algebras}

In this section we apply our formal algebraic uniformization result,
Corollary \ref{switching corollary}, to the Lie algebra $\mathfrak{g}$ 
consisting of certain formal infinite series with coefficients in a 
$\mathbb{Z}$-graded Lie algebra $\mathfrak{p}$.  We then give applications  
for $\mathfrak{p}$ an affine Lie algebra, the Virasoro algebra and a 
Grassmann envelope of the $N = 1$ Neveu-Schwarz algebra.  For the latter 
two applications, we 
discuss the importance of these results to conformal and superconformal 
field theory, and we point out that the result also applies to Grassmann
envelopes of the $N > 1$ Neveu-Schwarz algebras.

We continue to work over our field $\mathbb{F}$ of characteristic
zero.  Let $\mathfrak{p} = \coprod_{j \in \mathbb{Z}} \mathfrak{p}_j$
be a $\mathbb{Z}$-graded Lie algebra. (Note that if, for
example, $\mathfrak{p}$ is given as $\frac{1}{T} \mathbb{Z}$-graded
for some positive integer $T$, then by regrading, we can always
consider $\mathfrak{p}$ as $\mathbb{Z}$-graded.)  We consider the
following vector space decomposition of $\mathfrak{p}$:
\[\mathfrak{p} = \left( \coprod_{j < 0} \mathfrak{p}_j \right) \oplus
\mathfrak{p}_0\oplus 
\left( \coprod_{j > 0} \mathfrak{p}_j \right) = \mathfrak{p}^-
\oplus \mathfrak{p}_0\oplus \mathfrak{p}^+ .\] 

Let ${\cal A}_j$ (resp., ${\cal B}_j$) be commuting formal variables
for $j \in \mathbb{Z}$, and $j > 0$ (resp., $j < 0$).  Consider the
corresponding polynomial algebra $\mathbb{F}[{\cal A}, {\cal B}]$.  We
define the {\it order} of each ${\cal A}_j$ and each ${\cal B}_j$ to 
be one. This induces three gradings by nonnegative integers, called 
the {\it order in the ${\cal A}_j$'s}, the {\it order in the 
${\cal B}_j$'s} and the {\it total order in the ${\cal A}_j$'s and 
${\cal B}_j$'s}, defined in the obvious ways.

Now consider the space $\mathfrak{p}[{\cal A}, {\cal B}]$ of
polynomials in the ${\cal A}_j$'s and ${\cal B}_j$'s with coefficients
in $\mathfrak{p}$, equipped with the three gradings by order.  Also
consider the corresponding space $\mathfrak{p}[[{\cal A}, {\cal B}]]$
of formal power series.  Make $\mathfrak{p}[{\cal A},{\cal B}]$ and
$\mathfrak{p}[[{\cal A},{\cal B}]]$ into Lie algebras in the canonical
ways, and note that $\mathfrak{p}[[{\cal A},{\cal B}]]$ has a vector
space decomposition 
\begin{equation}\label{decomp}
\mathfrak{p}[[{\cal A},{\cal B}]]=\mathfrak{p}^-[[{\cal A}, 
{\cal B}]]\oplus \mathfrak{p}_{0}[[{\cal A}, 
{\cal B}]]\oplus \mathfrak{p}^+[[{\cal A}, 
{\cal B}]]
\end{equation}
into three Lie subalgebras $\mathfrak{p}^-[[{\cal A}, 
{\cal B}]]$, $\mathfrak{p}_{0}[[{\cal A}, 
{\cal B}]]$ and $\mathfrak{p}^+[[{\cal A}, 
{\cal B}]]$.

We are now ready to apply our results {}from the preceding section.  We
want to exponentiate suitable elements of the form $(p_1, p_2,...) \in
\prod_{j>0} \mathfrak{p}_j$, with possibly infinitely many $p_j$'s
nonzero, and corresponding elements of $\prod_{j<0} \mathfrak{p}_j$.
But products of such exponentials are not well defined even if we use
formal variables such as $s$ and $t$ as in the previous section.
However, if we use an infinite number of formal variables, one for
each homogeneous subspace, then we can multiply the exponentials.

Therefore, fix 
\begin{equation}
(p_1,p_2,...) \in \prod_{j>0} \mathfrak{p}_j \qquad  \mathrm{and} \qquad
(p_{-1},p_{-2},...) \in \prod_{j<0} \mathfrak{p}_j,
\end{equation}
and define corresponding elements 
\begin{equation}
g^+ = \sum_{j > 0} {\cal A}_j p_j \qquad \mathrm{and} \qquad g^- =
\sum_{j < 0} {\cal B}_j p_j
\end{equation}
of $\mathfrak{p}[[{\cal A},{\cal B}]]$.  We will consider exponentials
of these elements.

We will apply Corollary \ref{switching corollary} 
and Corollary \ref{switchingwordscorollary} to the Lie algebra 
$\mathfrak{p}[[{\cal A},{\cal B}]]$ with the decomposition 
$\mathfrak{p}[[{\cal A}, {\cal B}]] =\mathfrak{p}^- [[{\cal A}, 
{\cal B}]] \oplus (\mathfrak{p}_0[[{\cal A}, {\cal B}]] \oplus 
\mathfrak{p}^+[[{\cal A}, {\cal B}]])$ with corresponding projectins
$\pi^-: \mathfrak{p}[[{\cal A}, {\cal B}]] \rightarrow \mathfrak{p}^-
[[{\cal A}, {\cal B}]]$, and $\pi^{0,+}: \mathfrak{p}[[{\cal A}, {\cal B}]]
\rightarrow (\mathfrak{p}_0[[{\cal A}, {\cal B}]] \oplus \mathfrak{p}^+
[[{\cal A}, {\cal B}]])$, and then we will apply Theorem \ref{generalwordsthm}
to the Lie algebra $\mathfrak{p}_0
[[{\cal A}, {\cal B}]] \oplus \mathfrak{p}^+ [[{\cal A}, {\cal B}]]$ with 
the indicated decomposition and corresponding projections $\pi^0: 
\mathfrak{p}_0 [[{\cal A}, {\cal B}]] \oplus \mathfrak{p}^+ [[{\cal A}, 
{\cal B}]] \rightarrow \mathfrak{p}_0 [[{\cal A}, {\cal B}]]$ and $\pi^+: 
\mathfrak{p}_0 [[{\cal A}, {\cal B}]] \oplus \mathfrak{p}^+ [[{\cal A}, 
{\cal B}]] \rightarrow \mathfrak{p}^+ [[{\cal A}, {\cal B}]]$.

\begin{cor}\label{application corollary}
With the notation as above, there exist unique elements $\Psi^- 
\in \mathfrak{p}^-[[{\cal A}, {\cal B}]]$, $\Psi^{0} \in
\mathfrak{p}_{0}[[{\cal A}, {\cal B}]]$ and $\Psi^{+} \in 
\mathfrak{p}^{+}[[{\cal A}, {\cal B}]]$ of the form
\begin{eqnarray}
\Psi^- &=& \sum_{j < 0} {\cal B}_j p_j + {\cal Q}^-({\cal A}, {\cal
B}), \label{positive condition}\\ 
\Psi^{0} &=& {\cal Q}^0({\cal A}, {\cal B}),\label{zero condition}\\
\Psi^{+} &=& \sum_{j > 0} {\cal
A}_j p_j + {\cal Q}^+({\cal A}, {\cal B}), \label{negative
condition}
\end{eqnarray} 
where ${\cal Q}^- ({\cal A}, {\cal B}) \in \mathfrak{p}^-[[{\cal A}, 
{\cal B}]]$, ${\cal Q}^{0} ({\cal A}, {\cal B}) \in \mathfrak{p}_{0}
[[{\cal A}, {\cal B}]]$ and ${\cal Q}^{+} ({\cal A}, {\cal B}) \in 
\mathfrak{p}^{+}[[{\cal A}, {\cal B}]]$, and these elements contain 
only terms of order at least one in the ${\cal A}_j$'s and order at
least one in the ${\cal B}_j$'s, such that 
\begin{equation}\label{main-appl}
e^{g^+} e^{g^-} = e^{\Psi^-} e^{\Psi^{+}} e^{\Psi^{0}} 
\end{equation}
in $U(\mathfrak{p})[[{\cal A}, 
{\cal B}]]$. (Note that the right-hand side of (\ref{main-appl}) is
well defined for any elements of the form 
(\ref{positive condition})--(\ref{negative condition}).)
Moreover, 
\begin{equation}\label{wordsassertion}
\Psi^-, \Psi^0, \Psi^+ \in 
\mathbb{F} \langle \langle g^-, g^+; \pi^-, \pi^0, \pi^+  \rangle \rangle
\end{equation}
(using obvious notation), and we have
\begin{eqnarray}
\Psi^- &=& \sum_{j < 0} {\cal B}_j p_j + \sum_{\stackrel{j > 0, \; m <
0}{\mbox{\tiny $j + m < 0$}}} {\cal A}_j{\cal B}_m \left[ p_j, p_m\right]  
+ {\cal
P}^-({\cal A}, {\cal B}), \label{psi- condition}\\  
\Psi^{0} &=& \sum_{j > 0} {\cal A}_j{\cal B}_{-j} \left[p_j, p_{-j} 
\right] + {\cal P}^0({\cal A}, {\cal B}), \label{psi+,+ condition}\\ 
\Psi^{+} &=& \sum_{j > 0} {\cal A}_j p_j + \sum_{\stackrel{j > 0, \; m
< 0 }{\mbox{\tiny $j + m > 0$}}} 
{\cal A}_j{\cal B}_m \left[ p_j, p_m \right] +
{\cal P}^+({\cal A}, {\cal B}), \label{psi+,- condition}
\end{eqnarray} 
where ${\cal P}^-({\cal A}, {\cal B})$, 
${\cal P}^0({\cal A}, {\cal B})$,
${\cal P}^+({\cal A},
{\cal B})\in \mathfrak{p}[[{\cal
A}, {\cal B}]]$ each contain only terms of total order three or more
in the ${\cal A}_j$'s and ${\cal B}_j$'s, with order at least one in
the ${\cal A}_j$'s and at least one in the ${\cal B}_j$'s.
\end{cor}

{\it Proof:} \hspace{.2cm} The uniqueness of $\Psi^-$, $\Psi^{0}$ and
$\Psi^{+}$ is immediate {}from the argument in Remark
\ref{directuniqueness}, applied first to the decomposition
$$U(\mathfrak{p}) = U(\mathfrak{p}^-) \otimes U(\mathfrak{p}^+) \otimes
U(\mathfrak{p}_0)$$
and the corresponding projection
$$U(\mathfrak{p})[[{\cal A},{\cal B}]] \longrightarrow 
U(\mathfrak{p}^-)U(\mathfrak{p}^+)[[{\cal A},{\cal B}]].$$
This gives the uniqueness of the product $e^{\Psi^-} e^{\Psi^{+}}$,
and the analogous consideration of the decomposition
$$U(\mathfrak{p}^-)U(\mathfrak{p}^+) = U(\mathfrak{p}^-) \oplus 
U(\mathfrak{p}^-)U(\mathfrak{p}^+)\mathfrak{p}^+$$
inside $U(\mathfrak{p})$ then completes the uniqueness.

Note that $g^-\in \mathfrak{p}^-[[{\cal A}, {\cal B}]]$ and 
$g^+\in \mathfrak{p}^+[[{\cal A}, {\cal B}]]\subset \mathfrak{p}_0
[[{\cal A}, {\cal B}]] \oplus \mathfrak{p}^+[[{\cal A}, {\cal B}]]$.
By Corollaries \ref{switching corollary} and \ref{switchingwordscorollary},
there exist unique elements 
$\tilde{\Psi}^- \in (\mathfrak{p}^- [[{\cal A}, {\cal B}]])[[s,t]]$ and
$\tilde{\Psi}^{0, +} \in (\mathfrak{p}_{0}[[{\cal A}, {\cal B}]] \oplus 
\mathfrak{p}^{+}[[{\cal A}, {\cal B}]])[[s,t]]$  such that 
\[e^{tg^+} e^{sg^-} = e^{s\tilde{\Psi}^-} e^{t \tilde{\Psi}^{0, +}}, \]
and $s\tilde{\Psi}^-, t\tilde{\Psi}^{0,+} \in 
\mathbb{F} \langle \langle sg^-, tg^+; \pi^-, \pi^{0,+} \rangle 
\rangle \subseteq \mathbb{F} \langle \langle sg^-, tg^+; \pi^-, 
\pi^0, \pi^+ \rangle \rangle$, and these elements satisfy 
(\ref{lowest order G -}) and (\ref{lowest order G +}), respectively, 
applied to this case.  For any element of $\mathbb{F}\langle 
\langle sg^-, tg^+;\pi^-,\pi^{0,+}\rangle\rangle$, 
we can substitute $1$ for both $s$ and 
$t$ and the result is an element of $\mathfrak{p}[[\mathcal{A}, 
\mathcal{B}]]$. Thus $\Psi^- = s \tilde{\Psi}^- \left.\right|_{s = t
= 1}$ and $\Psi^{0,+} = t \tilde{\Psi}^{0, +} \left.\right|_{s = t =
1}$ are well-defined elements of $\mathfrak{p}^-[[{\cal A},{\cal B}]]$ 
and $\mathfrak{p}_0 [[{\cal A},{\cal B}]] \oplus \mathfrak{p}^+
[[{\cal A},{\cal B}]] $, respectively, and $e^{\Psi^-}$
and $e^{\Psi^{0,+}}$ are well defined in $U(\mathfrak{p})[[{\cal A},
{\cal B}]]$.  Moreover, $\Psi^-$ and $\Psi^{0,+}$ satisfy 
the analogues of (\ref{lowest 
order G -}) and (\ref{lowest order G +}), respectively, applied to 
this case, and thus $\Psi^-$ satisfies (\ref{psi- condition}).

Now let $h^+ = \pi^{+} (\Psi^{0,+})$ and $h^0 = \pi^{0} (\Psi^{0,+})$.
By Theorems \ref{main theorem} and \ref{generalwordsthm}, 
there exist unique $\tilde{\Psi}^+ \in 
(\mathfrak{p}^+ [[{\cal A},{\cal B}]])[[s,t]]$ and $\tilde{\Psi}^0 \in 
(\mathfrak{p}_0 [[{\cal A},{\cal B}]]) [[s,t]]$ such that
\[e^{s\tilde{\Psi}^+} e^{t\tilde{\Psi}^0} = e^{sh^+ + th^0},\]
and $s\tilde{\Psi}^+, t\tilde{\Psi}^0 \in \mathbb{F}\langle 
\langle sh^+,th^0; \pi^+, \pi^0 \rangle \rangle$ and $s\tilde{\Psi}^+$ 
and $t\tilde{\Psi}^0$ satisfy (\ref{lowest order F -}) and 
(\ref{lowest order F +}), respectively, applied to this case.  
Again we can substitute 
$1$ for both $s$ and $t$ and the result is an element of $\mathfrak{p}
[[\mathcal{A}, \mathcal{B}]]$. Thus we have that $\Psi^+ = s 
\tilde{\Psi}^+ \left.\right|_{s = t = 1}$ and $\Psi^0 = t \tilde{\Psi}^0 
\left.\right|_{s = t = 1}$ are well defined elements of 
$\mathfrak{p}^+[[{\cal A},{\cal B}]]$ and $\mathfrak{p}_0 [[{\cal A},
{\cal B}]]$, respectively, and $e^{\Psi^-}$ and $e^{\Psi^{0,+}}$ are 
well defined in $U(\mathfrak{p})[[{\cal A}, {\cal B}]]$.  Moreover,
$\Psi^+, \Psi^0 
\in \mathbb{F}\langle \langle h^+, h^0; \pi^+, \pi^0 \rangle 
\rangle \subseteq \mathbb{F}\langle \langle g^-, g^+; \pi^-,
\pi^+, \pi^0 \rangle \rangle$, and
by (\ref{lowest order F -}) and (\ref{lowest order F +}) applied to this
case, $\Psi^0$ and $\Psi^+$ satisfy (\ref{psi+,+ condition}) and 
(\ref{psi+,- condition}), respectively. 
\pfbox

Now we specialize to the following situations:

\begin{app}\label{affine} {Affine Lie algebras.} 
{\em Let $\mathfrak{l}$ be a finite dimensional Lie algebra equipped with
an  $\mathfrak{l}$-invariant symmetric bilinear form $(\cdot, \cdot)$
and consider the corresponding affine Lie algebra
$\hat{\mathfrak{l}} = \mathfrak{l} \otimes \mathbb{C}[x, x^{-1}]
\oplus \mathbb{C}c,$  with
commutation relations
\[ [g \otimes x^m, h \otimes x^n] = [g,h] \otimes x^{m + n} + 
(g, h) m \delta_{m + n,0}{\bf k} , \]
and ${\bf k}$ central. Consider also the natural $\mathbb{Z}$-grading 
$\hat{\mathfrak{l}} = \coprod_{n \in \mathbb{Z}} \hat{\mathfrak{l}}_n$.

Fix $h_j \in \mathfrak{l}$ for $j \in \mathbb{Z} \smallsetminus \{0\}$
and let ${\cal A}_j$ for $j>0$ and ${\cal B}_j$ for $j<0$ be commuting
formal variables.  We can now apply Corollary \ref{application
corollary} with $\mathfrak{p} = \hat{\mathfrak{l}}$,
$\mathfrak{p}_j = \hat{\mathfrak{l}}_j$ for $j \in
\mathbb{Z}$ and $p_j = h_j \otimes x^j$ for $j \in \mathbb{Z}
\smallsetminus \{0\}$.  Thus
\[ g^+ = \sum_{j > 0} {\cal A}_j \; h_j \otimes x^j \qquad \mathrm{and} \qquad
g^- = \sum_{j < 0} {\cal B}_j \; h_j \otimes x^j, \]
and we have that 
there exist unique $\Psi^- \in \hat{\mathfrak{l}}^-[[{\cal A}, 
{\cal B}]]$, $\Psi^{0} \in \hat{\mathfrak{l}}_0[[{\cal A}, {\cal B}]]$, 
and $\Psi^{+} \in \hat{\mathfrak{l}}^+[[{\cal A}, {\cal B}]]$ 
satisfying (\ref{positive condition}) - (\ref{negative
condition}) such that
\begin{equation}\label{affine equation}
\exp \left( \sum_{j > 0} {\cal A}_j \; h_j \otimes x^j \right) \exp
\left( \sum_{j < 0} {\cal B}_j \; h_j \otimes x^j \right) = e^{\Psi^-}
e^{\Psi^{+}} e^{\Psi^{0}}.
\end{equation}
Also, formula (\ref{wordsassertion}) holds.
Furthermore, we can write $\Psi^{0} \in \hat{\mathfrak{l}}_0 [[{\cal
A}, {\cal B}]]$ as $\Psi^{0} = \Psi^{0}_{\mathfrak{l}} + \Psi^{0}_{\bf k}$ 
where
$\Psi^{0}_{\mathfrak{l}}\in \mathfrak{l}[[{\cal A}, {\cal B}]]$ and
$\Psi^{0}_{\bf k}\in (\mathbb{C}{\bf k})[[{\cal A}, {\cal B}]]$.  
But since ${\bf k}$
is central, we can write the last exponential in equation (\ref{affine
equation}) as
\[ e^{\Psi^{0}} = \exp \left(\Psi^{0}_{\mathfrak{l}}+ \Psi^{0}_{\bf k}
\right) = 
\exp \left(\Psi^{0}_{\mathfrak{l}}\right) \exp \left(\Psi^{0}_{\bf k}
\right). \]

In addition, we have
\begin{eqnarray}
\Psi^- &=& \sum_{j < 0} {\cal B}_j \; h_j \otimes x^j + \hspace{3in}
\nonumber \\ 
& & \hspace{.8in} \sum_{\stackrel{j > 0, \; m < 0 }{\mbox{\tiny
$j + m <0$}}} {\cal
A}_j {\cal B}_m \left[h_j, h_m \right] \otimes x^{j + m} + {\cal
P}^-({\cal A}, {\cal B}), \\ 
\Psi^{+} &=& \sum_{j > 0} {\cal A}_j \; h_j \otimes x^j + \nonumber
\\ 
& & \hspace{.8in} \sum_{\stackrel{j > 0, \; m < 0 }{\mbox{\tiny
$j + m > 0$}}} {\cal
A}_j {\cal B}_m \left[ h_j, h_m \right] \otimes x^{j + m} + {\cal
P}^+({\cal A}, {\cal B}), \\ 
\Psi^{0}_{\mathfrak{l}}
&=& \sum_{j > 0} {\cal A}_j {\cal B}_{-j} \left[ h_j,
h_{-j} \right] + {\cal P}^{0}_{\mathfrak{l}}({\cal A}, {\cal B}), \\ 
\Psi^{0}_{\bf k}&=& \sum_{j > 0} {\cal A}_j {\cal B}_{-j} \left( h_j,
h_{-j} \right) j{\bf k} + {\cal P}^{0}_{\bf k} ({\cal A}, {\cal B})
\hspace{.6in}
\end{eqnarray} 
where ${\cal P}^-({\cal A}, {\cal B})$, ${\cal P}^+({\cal A}, {\cal B})$,
${\cal P}^{0}_{\mathfrak{l}}({\cal A}, {\cal B})$, 
${\cal P}^{0}_{\bf k}({\cal A}, {\cal B})
\in \hat{\mathfrak{l}}[[{\cal A}, {\cal B}]]$ each contain only terms
of total order three or more in the ${\cal A}_j$'s and ${\cal B}_j$'s,
with order at least one in the ${\cal A}_j$'s and at least one in the
${\cal B}_j$'s.  }
\end{app}

\begin{app}\label{virasoro example} The Virasoro algebra.  
{\em Take $\mathfrak{p} = \mathfrak{v}$ to be the Virasoro algebra.
With the usual basis, $L_{m}$ for $m\in \mathbb{Z}$ and $c$,
the commutation relations for
$\mathfrak{v}$ are
\begin{eqnarray*}
\left[ L_m , L_n\right] &=& (m - n) L_{m + n} + \frac{1}{12}(m^3 - m) 
\delta_{m + n, 0} c ,\\
\left[ L_m, c \right] &=& 0 
\end{eqnarray*}
for $m,n \in \mathbb{Z}$. Consider the natural $\mathbb{Z}$-grading 
$\mathfrak{v}=\coprod_{j\in \mathbb{Z}}\mathfrak{v}_{j}$, where 
$\mathfrak{v}_{j}=\mathbb{C}L_{j}$ for $j\in \mathbb{Z}\setminus \{0\}$ and 
$\mathfrak{v}_{0}=\mathbb{C}L_{0}\oplus \mathbb{C}c$. 

Now take $p_j = L_j$ for $j \in \mathbb{Z} \smallsetminus \{0\}$, 
and as usual 
let
${\cal A}_j$ and ${\cal B}_{-j}$ be commuting formal variables for $j
\in \mathbb{Z}$, $j > 0$.  (Note that we could take $p_j = c_j L_j$ where
$c_j$ are complex variables.  However, in expressions such as
${\cal A}_j c_j L_j$, we can always absorb the complex variables $c_j$
into the formal variables ${\cal A}_j$.  In fact, in applications such
as in \cite{H1} - \cite{H3}, under suitable conditions, one eventually
wants to specialize the formal variables to be complex numbers.  Thus we
have not sacrificed any generality by setting $c_j = 1$.)  

Now take
\[g^+ = \sum_{j > 0} {\cal A}_j L_j \qquad \mathrm{and} \qquad 
g^- = \sum_{j < 0} {\cal B}_j L_j.\]  
Applying Corollary \ref{application corollary}, we see that there 
exist unique 
$\Psi^- \in \mathfrak{v}^{-}[[{\cal A}, {\cal B}]]$, $\Psi^{0} \in 
\mathfrak{v}_{0}[[{\cal A}, {\cal B}]]$ and $\Psi^{+} \in
\mathfrak{v}^{+}[[{\cal A}, {\cal B}]]$ satisfying (\ref{positive condition}) 
- (\ref{negative condition}), such that
\begin{equation}\label{virasoro}
\exp \left( \sum_{j > 0} {\cal A}_j L_j \right) \exp \left( \sum_{j <
0} {\cal B}_j L_j \right) = e^{\Psi^-} e^{\Psi^{+}} e^{\Psi^{0}},
\end{equation}
and formula (\ref{wordsassertion}) holds.
Let us write
\[\Psi^- = \sum_{j<0} \Psi_j L_j, \qquad \Psi^{+} = \sum_{j>0} \Psi_j L_j\]
and 
\[\Psi^{0} = \Psi_{0} L_0 + \Gamma c\] 
where $\Psi_j, 
\Gamma \in \mathbb{C}[[{\cal A}, {\cal B}]]$ for
$j\in \mathbb{Z}$.
Since $c$ is central, (\ref{virasoro}) is equal to
\[ \exp \left( \sum_{j < 0} \Psi_j L_j \right) \exp \left(\sum_{j > 0}
\Psi_{j} L_j \right) e^{\Psi_{0}L_0} e^{\Gamma c}, \]
and for $j > 0$, we have 
\begin{eqnarray}
\Psi_{-j} &=& {\cal B}_{-j} + \sum_{m > j} {\cal A}_{-j + m}
{\cal B}_{-m} (-j + 2m) + {\cal P}_{-j} ({\cal A}, {\cal B}),
\label{v1} \\ 
\Psi_j &=& {\cal A}_j + \sum_{m > 0} {\cal A}_{j + m} 
{\cal B}_{-m} (j + 2m) + {\cal P}_j ({\cal A}, {\cal B}), \label{v2}
\\ 
\Psi_{0} &=& \sum_{m > 0} 2 {\cal A}_m {\cal B}_{-m} m + 
{\cal P}_{0} ({\cal A}, {\cal B}), \label{v3}\\
\Gamma &=& \sum_{m > 0} {\cal A}_m {\cal B}_{-m} 
\frac{(m^3 - m)}{12} + \Gamma_{0}({\cal A}, {\cal B}), \label{v4}
\end{eqnarray}
where ${\cal P}_{j} ({\cal A}, {\cal B}), \Gamma_{0}({\cal
A}, {\cal B}) \in \mathbb{C}[[{\cal A}, {\cal B}]]$, for $j\in
\mathbb{Z}$,
contain only terms of total order three or more in the ${\cal A}_m$'s
and ${\cal B}_m$'s with order at least one in the ${\cal A}_m$'s and at
least one in the ${\cal B}_m$'s.

\begin{rema}\label{uniformization v} 
{\em In conformal field theory, equation (\ref{virasoro}) corresponds
to calculating the uniformizing function to obtain a canonical sphere
with tubes {}from the sewing together of two canonical spheres with
tubes in the moduli space of spheres with tubes under global conformal
equivalence.  This moduli space along with the sewing operation is the
geometric structure underlying a geometric vertex operator algebra
\cite{H1} - \cite{H3}.  Equation (\ref{virasoro}) also corresponds to
a certain change of variables and ``normal ordering'' of the operators
$L_j$ generated by the Virasoro element in an (algebraic) vertex
operator algebra, where by ``normal ordering'' we mean ordering the
operators $L_j$ so as to first act by the operators $L_j$ for $j
> 0$ and then act by the operators $L_j$ for $j < 0$.  The
correspondence between these two procedures, one geometric and the other
algebraic, is necessary for the proof of the isomorphism between the
category of vertex operator algebras and the category of
geometric vertex operator algebras \cite{H2}.}
\end{rema}

\begin{rema}\label{necessary data v} 
{\em The results about the formal series $\Psi_j, \Gamma\in
\mathbb{C}[[{\cal A}, {\cal B}]]$ for $j \in \mathbb{Z}$, given in 
equations (\ref{v1}) - (\ref{v4})
above---the explicit results about the lowest order terms and the
qualitative information about the higher order terms---are exactly the
results necessary for the proof of the isomorphism between the
category of geometric vertex operator algebras and the category of
algebraic vertex operator algebras.  Equivalent results were
proved by Huang, first in \cite{H1} and then in Theorem 2.2.4\footnote{We
take this opportunity to correct a misprint in the formulas 
(2.2.11) and (2.2.12) of Theorem 
2.2.4 in \cite{H3}.  The first two terms in the right-hand side of (2.2.11) 
should be replaced by $\alpha_0^{-1}(f_{\cal B}^{(2)})^{-1}(\frac{1}
{\alpha_0 x})$ 
and the first two terms in the right-hand side of (2.2.12) should be replaced
by $(f_{{\cal A},\alpha_0}^{(1)})^{-1} (x)$.}, Proposition
2.2.5\footnote
{There 
is a misprint in equation (2.2.27) of Proposition 2.2.5 in \cite{H3}.  The 
first term in the right-hand side of (2.2.27) should be $-\alpha_0^{-j} 
{\cal B}_j$.}, Proposition 4.2.1 and Corollary 4.2.2 of \cite{H3}.  However, 
the quantitative information in equations (\ref{v1}) - (\ref{v4}) is much more 
explicit than the equivalent information given in
\cite{H3}. The main difference between equations (\ref{v1}) -
(\ref{v4}) and the analogous results given in \cite{H3} is that 
in (\ref{v1}) - (\ref{v4}), the terms of total order two in the
${\cal A}_j$'s and ${\cal B}_j$'s are given 
explicitly while in \cite{H3} this information for the $\Psi_j$'s and $\Gamma$
is presented in (the corrected forms of) equations (2.2.11) and
(2.2.12) (see footnote 4) using
the representation of the Virasoro algebra given by $L_n = - x^{n + 1}
\frac{d}{d x} \in \mathrm{End} (\mathbb{C}[x, x^{-1}])$ and $c = 0$.
In order to recover equations 
(\ref{v1}) - (\ref{v4}) above {}from the 
results in \cite{H3}, one must perform several operations, pick out 
coefficients, and then
use Proposition 4.2.1 and Corollary 4.2.2 in \cite{H3}, allowing one
to lift the results {}from the particular representation to the algebra.  For
example, a shorter and more straightforward proof of Proposition 3.5.2
in \cite{H3} than that originally given can be obtained using
equations (\ref{v1}) - (\ref{v4}) above.  This proposition states that
the meromorphic tangent space of the moduli space of spheres with one
incoming tube and one outgoing tube carries the structure of a
Virasoro algebra with central charge zero. }
\end{rema}}
\end{app}

We now show how Corollary \ref{application corollary} can be
applied to Lie superalgebras.  We will use the notion of a Grassmann
envelope of a Lie superalgebra. These are  Lie algebras 
to which we will apply Corollary \ref{application corollary}.

We will be interested in $\frac{1}{2}\mathbb{Z}$-graded vector spaces
also equipped with a compatible $\mathbb{Z}_2$-grading.  To
distinguish, we will denote the $\mathbb{Z}_2$-grading using
superscripts.  For a $\mathbb{Z}_2$-graded vector space $V
= V^0 \oplus V^1$, define the {\it sign function} $\eta$ on the
homogeneous subspaces of $V$ by $\eta(v) = i$ for $v \in V^i$, $i =
0,1$.  If $\eta(v) = 0$, we say that $v$ is {\it even}, and if
$\eta(v) = 1$, we say that $v$ is {\it odd}.

A {\it superalgebra} is an (associative) algebra $A$ (with identity $1
\in A$), such that
\begin{eqnarray*}
&(i)& \mbox{$A$ is a $\mathbb{Z}_2$-graded algebra}\\ &(ii)& \mbox{$ab
= (-1)^{\eta(a) \eta(b)} ba$ for $a,b$ homogeneous in $A$.}
\hspace{1.5in}
\end{eqnarray*}
For example, the exterior (or Grassmann) 
algebra $\bigwedge(V)$ over a vector space $V$ is naturally 
a superalgebra.

A $\mathbb{Z}_2$-graded vector space $\mathfrak{q}$ is said to be a
{\it Lie superalgebra} if it has a bilinear operation $[\cdot,\cdot]$
such that for $u,v$ homogeneous in $\mathfrak{q}$,
\begin{eqnarray*}
&(i)& [u,v] \in {\mathfrak q}^{(\eta(u) + \eta(v)) \mathrm{mod} \; 2}\\ 
&(ii)& [u,v] = -(-1)^{\eta(u)\eta(v)}[v,u] \hspace{1.75in} 
\mbox{(skew-symmetry)}\\
&(iii)& (-1)^{\eta(u)\eta(w)}[[u,v],w] + (-1)^{\eta(v)\eta(u)}[[v,w],u] \\
& & \hspace{1.35in} + \; (-1)^{\eta(w)\eta(v)}[[w,u],v] = 0. \qquad 
\mbox{(Jacobi identity)}
\end{eqnarray*}

\begin{rema}\label{envelope}
{\em Given a Lie superalgebra $\mathfrak{q}$ and a superalgebra $A$, $(A^0 
\otimes \mathfrak{q}^0) \oplus (A^1 \otimes
\mathfrak{q}^1)$ is a Lie algebra with bracket given by
\begin{equation}\label{super v. non} 
[au , bv] = (-1)^{\eta(b) \eta(u)} ab [u,v]
\end{equation}
(with obvious notation), where we have suppressed the tensor product
symbol.  Note that the bracket on the left-hand side of (\ref{super
v. non}) is a Lie algebra bracket, and the bracket on the right-hand
side is a Lie superalgebra bracket.  If $A = \bigwedge(V)$ for some
vector space $V$, then this Lie
algebra is called the {\it Grassmann envelope of the Lie superalgebra
$\mathfrak{q}$ associated with $A$}.}
\end{rema}

Consider a Lie superalgebra $\mathfrak{q}$ that also has a compatible
$\mathbb{Z}$-grading. Fix a Grassmann algebra $A=\bigwedge(V)$.  We can
now apply Corollary \ref{application corollary} to the Grassmann
envelope of $\mathfrak{q}$ associated with $A$. In addition, we can
apply  Corollary \ref{application corollary} to the Grassmann
envelope of a $\frac{1}{T} \mathbb{Z}$-graded Lie superalgebra
associated with $A$ by regrading.

\begin{app}\label{N = 1} The $N = 1$ Neveu-Schwarz algebra.
{\em The $N=1$ Neveu-Schwarz Lie superalgebra, $\mathfrak{ns}$, is a
superextension of the Virasoro algebra.  Thus $\mathfrak{ns}^0$ is the
Virasoro algebra $\mathfrak{v}$ as in Application \ref{virasoro
example}, and for $j \in \mathbb{Z} + \frac{1}{2}$, we have
$\mathfrak{ns}_j = \mathbb{C} G_j$.  The remaining supercommutation
relations are
\begin{eqnarray*} 
\left[ G_{m + \frac{1}{2}}, L_n \right] &=& \left( m - \frac{n - 1}{2} 
\right) G_{m + n + \frac{1}{2}} \\
\left[ G_{m + \frac{1}{2}}, G_{n - \frac{1}{2}} \right] &=& 2L_{m + n} + 
\frac{1}{3} (m^2 + m) \delta_{m + n, 0} c \\
\left[ G_{m + \frac{1}{2}}, c \right] &=& 0 
\end{eqnarray*}
for $m,n \in \mathbb{Z}$.

Take $\mathfrak{p} = (A^0 \otimes \mathfrak{ns}^0) \oplus
(A^1 \otimes \mathfrak{ns}^1)$, $p_j = L_j$ for $j \in
\mathbb{Z} \smallsetminus \{0\}$, and  $p_{j - \frac{1}{2}} = a_{j
- \frac{1}{2}} G_{j - \frac{1}{2}}$ for 
 $j \in \mathbb{Z}$ where $a_{j - \frac{1}{2}} \in
A^1$.    Let ${\cal A}_j$ and ${\cal B}_{-j}$ be 
commuting formal variables for $j \in \frac{1}{2} \mathbb{Z}$, 
$j>0$, and set
\begin{eqnarray*}
g^+ &=& \sum_{j \in \mathbb{Z}_+} \left( {\cal A}_j L_j + {\cal A}_{j -
\frac{1}{2}} a_{j - \frac{1}{2}} G_{j - \frac{1}{2}} \right) \\
g^- &=& \sum_{j \in -\mathbb{Z}_+} \left( {\cal B}_j
L_j + {\cal B}_{j + \frac{1}{2}} a_{j + \frac{1}{2}} G_{j +
\frac{1}{2}} \right)
\end{eqnarray*} 
($\mathbb{Z}_+$ denoting the positive integers).  By Corollary
\ref{application corollary}, there exist unique $\Psi^- \in \mathfrak{p}^{-}
[[{\cal A}, {\cal B}]]$, $\Psi^{0} \in \mathfrak{p}_{0}[[{\cal A}, 
{\cal B}]]$ and $\Psi^{+} \in \mathfrak{p}^{+}[[{\cal A}, {\cal B}]]$, 
satisfying (\ref{positive condition}) -
(\ref{negative condition}) such that
\begin{equation}\label{neveu-schwarz}
e^{g^+} e^{g^-} = e^{\Psi^-} e^{\Psi^{+}} e^{\Psi^{0}},
\end{equation}    
and formula (\ref{wordsassertion}) holds.
Since $\{L_j, G_{j- \frac{1}{2}}\}_{j \in \mathbb{Z}}
\cup \{c\}$ is a basis for the Neveu-Schwarz algebra, we can write
\[\Psi^- = \sum_{j \in -\mathbb{Z}_+} \left(\Psi_j L_j + 
\Psi_{j+\frac{1}{2}}G_{j + \frac{1}{2}}
\right), \qquad \Psi^{+} = \sum_{j \in \mathbb{Z}_+} \left(\Psi_j
L_j + \Psi_{j - \frac{1}{2}}G_{j - \frac{1}{2}} \right)\] 
and
\[\Psi^{0} = \Psi_0 L_0 + \Gamma c \] 
where $\Psi_j, \Gamma \in A[[{\cal A}, {\cal B}]]$ for
$j\in \frac{1}{2}\mathbb{Z}$.  Since $c$ is central, 
\begin{equation}
e^{\Psi^{0}} = \exp \left( \Psi_0 L_0 + \Gamma c \right) =
\exp \left( \Psi_0 L_0 \right) \exp \left( \Gamma c \right),
\end{equation} 
and for $j \in \mathbb{Z}_+$, we have
\begin{eqnarray}
\Psi_{-j} &=& {\cal B}_{-j} + \sum_{m > j} \left( {\cal
A}_{-j + m} {\cal B}_{-m} (-j + 2m) \right. \label{ns1} \\ & &
\hspace{.5in} \left. - \; 2 {\cal A}_{-j + m - \frac{1}{2}} {\cal
B}_{- m + \frac{1}{2}} a_{-j + m - \frac{1}{2}} a_{- m + \frac{1}{2}}
\right) \nonumber + {\cal P}_{-j} ({\cal A}, {\cal B}), \nonumber \\
\Psi_{-j + \frac{1}{2}} &=& {\cal B}_{- j + \frac{1}{2}} a_{- j +
\frac{1}{2}} + \! \! \sum_{m > j} \! \left( {\cal A}_{-j +
m} {\cal B}_{- m + \frac{1}{2}} a_{- m + \frac{1}{2}}
\left(-\frac{j}{2}+\frac{3m}{2}-\frac{1}{2} \right) \right. \label{ns2} \\ & &
\left. \hspace{-.2in} + \; {\cal A}_{-j + m - \frac{1}{2}} {\cal
B}_{-m + 1} a_{-j + m - \frac{1}{2}} \left(-j + \frac{3m}{2} - 1
\right) \right) + {\cal P}_{-j + \frac{1}{2}} ({\cal A}, {\cal B}),
\hspace{.3in} \nonumber \\  
\Psi_j &=& {\cal A}_j + \sum_{m \in \mathbb{Z}_+} \left( {\cal A}_{j + m}
{\cal B}_{-m} (j + 2m) \right. \label{ns3}\\
& & \hspace{.8in} \left. - \; 2 {\cal A}_{j + m - \frac{1}{2}} 
{\cal B}_{-m + \frac{1}{2}} a_{j + m - \frac{1}{2}} 
a_{-m + \frac{1}{2}} \right) + {\cal P}_j ({\cal A}, {\cal B}), 
\nonumber \\
\Psi_{j - \frac{1}{2}} &=& {\cal A}_{j - \frac{1}{2}} a_{j -
\frac{1}{2}} + \sum_{m \in \mathbb{Z}_+} \left(
{\cal A}_{j + m - 1} {\cal B}_{- m + \frac{1}{2}} a_{- m + \frac{1}{2}}
\left(\frac{j}{2} + \frac{3m}{2} - 1 \right) \right. \label{ns4} \\ 
& & \left. \hspace{.3in} + \; {\cal A}_{j + m - \frac{1}{2}} {\cal
B}_{-m} a_{j + m - \frac{1}{2}} \left(j + \frac{3m}{2} - \frac{1}{2} 
\right) \right) + {\cal P}_{j - \frac{1}{2}} ({\cal A}, {\cal B}), 
\nonumber  \\
\Psi_0 &=& \sum_{m \in \mathbb{Z}_+} \left( 2{\cal A}_m {\cal B}_{-m} m
- 2 {\cal A}_{m - \frac{1}{2}} {\cal B}_{-m + \frac{1}{2}} a_{m -
\frac{1}{2}} a_{-m + \frac{1}{2}} \right) \label{ns5}\\
& & \hspace{3.2in} + \; {\cal P}_0 ({\cal A}, {\cal B}), \nonumber \\ 
\Gamma &=& \sum_{m \in \mathbb{Z}_+} \left( {\cal A}_m {\cal B}_{-m}
\frac{(m^3 - m)}{12} \right. \label{ns6}\\
& & \left. \hspace{.5in} - \; {\cal A}_{m - \frac{1}{2}} 
{\cal B}_{-m + \frac{1}{2}} a_{m - \frac{1}{2}} a_{-m + \frac{1}{2}} 
\frac{(m^2 - m)}{3} \right) + \Gamma_0 ({\cal A}, {\cal B}), \nonumber 
\end{eqnarray}
where ${\cal P}_l ({\cal A}, {\cal B})$, $\Gamma_0 ({\cal A},{\cal
B}) \in A[[{\cal A}, {\cal B}]]$, for $l \in \frac{1}{2}
\mathbb{Z}$, contain only terms of
total order three or more in the ${\cal A}_m$'s and ${\cal B}_m$'s
with order at least one in the ${\cal A}_m$'s and at least one in the
${\cal B}_m$'s, for $m \in \frac{1}{2}\mathbb{Z}$.

\begin{rema}\label{uniformization ns} 
{\em (cf. Remark \ref{uniformization v})
In $N = 1$ superconformal field theory, equation
(\ref{neveu-schwarz}) corresponds to calculating the uniformizing
function to obtain a canonical supersphere with tubes {}from the sewing
together of two canonical superspheres with tubes in the moduli space
of superspheres with tubes under global superconformal equivalence.
This moduli space along with the sewing operation is the
supergeometric structure underlying a supergeometric vertex operator
superalgebra \cite{B1}, \cite{B2}.  Equation (\ref{neveu-schwarz})
also corresponds to a certain change of variables and ``normal
ordering'' of the operators $L_j, G_{j - \frac{1}{2}}$, $j \in
\mathbb{Z}$, generated by the Neveu-Schwarz element in an algebraic $N
= 1$ vertex operator superalgebra, where by ``normal ordering'' we
mean ordering the operators $L_j, G_{j - \frac{1}{2}}$ so as to first
act by the operators $L_j, G_{j - \frac{1}{2}}$ for $j > 0$ and
then act by the operators $L_j, G_{j + \frac{1}{2}}$ for $j < 0$.
The correspondence between these two procedures, one geometric and the other
algebraic, is necessary for the proof of the isomorphism between the
category of $N = 1$ vertex operator superalgebras and the category of
$N = 1$ supergeometric vertex operator superalgebras \cite{B1}.}
\end{rema}

\begin{rema}\label{necessary data ns} 
{\em (cf. Remark \ref{necessary data v})
The results about the formal series $\Psi_j, \Gamma \in
A[[{\cal A}, {\cal B}]]$, for $j \in
\frac{1}{2} \mathbb{Z}$, given in equations (\ref{ns1}) - (\ref{ns6})
above---the explicit results about the lowest order terms and the
qualitative information about the higher order terms---are exactly the
results necessary for the proof of the isomorphism between the
category of $N = 1$ vertex operator superalgebras and the category of
$N = 1$ supergeometric vertex operator superalgebras \cite{B1}.
Equivalent results were proved by Barron in Theorem 2.3.4\footnote{There 
is a misprint in formulas (2.49) and (2.50) of Theorem 2.3.4 in \cite{B1}.  
In the lowest order terms for $\bar{F}^{(1)}$, the first two terms of the 
even part and the first term and last two terms of the odd part of the 
right-hand side of formula (2.49) should be replaced by $\left. \varphi 
\bar{F}^{(1)} (x,\varphi) \right|_{({\cal A}, {\cal M}) = \mathbf{0}}$ 
given in equation (2.45) of the theorem.  And in the lowest order terms 
for $\bar{F}^{(2)}$, the first two terms of the even part and the first 
three terms of the odd part of the right-hand side of formula (2.50) 
should be replaced by $\left.\varphi \bar{F}^{(2)} (x,\varphi)
\right|_{({\cal B}, {\cal N}) = \mathbf{0}}$ given in equation (2.48) of 
the theorem.}, Proposition 2.3.6, Proposition 2.6.1 and Corollary 2.6.2 
in \cite{B1}.  However, the quantitative information in equations (\ref{ns1}) 
- (\ref{ns6}) is in a slightly different form than the equivalent information 
given in \cite{B1} and \cite{B2}.  In \cite{B1} and \cite{B2} ``odd'' formal 
variables $M_{j - \frac{1}{2}}$ (resp., $N_{-j - \frac{1}{2}}$) are used 
instead of the composite expressions ${\cal A}_{j - \frac{1}{2}} a_{j -
\frac{1}{2}}$ (resp., ${\cal B}_{j + \frac{1}{2}} a_{j +
\frac{1}{2}}$) found in equations (\ref{ns1}) - (\ref{ns6}) and
consisting of an even formal variable and an odd Grassmann variable.
The ``odd'' formal variables $M_{j - \frac{1}{2}}$ (resp., $N_{-j +
\frac{1}{2}}$) carry the same information as the corresponding
composite expressions in the present work and are ``odd'' in the sense
that they anticommute with each other and odd elements of
$\mathfrak{ns}$ and commute with even formal variables and even
elements of $\mathfrak{ns}$.  After taking into consideration this
notational change, one can see that the quantitative information in
equations (\ref{ns1}) - (\ref{ns6}) is much more explicit than the
equivalent information given in \cite{B1}. The main difference between
equations (\ref{ns1}) - (\ref{ns6}) and the analogous results first
proved in \cite{B1} is that in (\ref{ns1}) - (\ref{ns6}), the terms of 
total order two in the ${\cal A}_j$'s and ${\cal B}_j$'s are given 
explicitly, while in \cite{B1} this information for the $\Psi_j$'s and
$\Gamma$ is presented in (the corrected forms of) equations (2.49)
and (2.50) (see footnote 6) using a representation of the $N = 1$ Neveu-Schwarz algebra
in terms of superderivations in $\mathrm{End} (\mathbb{C}[x,
x^{-1}][\varphi])$ with $c = 0$, where $x$ is a formal (commuting)
variable and $\varphi$ is a formal anticommuting variable (see
Proposition 2.4.1 in \cite{B1}).  In order to recover equations 
(\ref{ns1}) - (\ref{ns6}) above {}from the results in \cite{B1}, 
one must perform several
operations, pick out coefficients, and then use Proposition 2.3.4 and
Corollary 2.6.2 in \cite{B1}, allowing one to lift {}from the particular
representation to the algebra.  For example, a shorter and more
straightforward proof of Proposition 3.11.1 in \cite{B1} than that
originally given can be obtained using equations (\ref{ns1}) -
(\ref{ns6}) above.  This proposition states that the supermeromorphic
tangent space of the moduli space of superspheres with one incoming
tube and one outgoing tube carries the structure of an $N = 1$
Neveu-Schwarz algebra with central charge zero. }
\end{rema} }
\end{app}

\begin{app}\label{bigger N} The Neveu-Schwarz algebras for $N > 1$.
{\em For Grassmann envelopes of other superextensions of the Virasoro
algebra, such as the $N = 2$ Neveu-Schwarz algebra and Neveu-Schwarz
algebras for higher $N$, the results of Corollary \ref{application
corollary} similarly apply.  These results have significance for the
corresponding superconformal field theories and vertex operator
superalgebras.}
\end{app}


\begin{thebibliography}{99}

\bibitem[B1]{B1} K. Barron, {\it The supergeometric interpretation of
vertex operator superalgebras}, Ph.D. thesis, Rutgers University,
1996.

\bibitem[B2]{B2} K. Barron, {\it A supergeometric interpretation of
vertex operator superalgebras}, Int. Math. Res. Notices, 1996 No. 9,
Duke University Press, 409--430.

\bibitem[H1]{H1} Y.-Z. Huang, {\it On the geometric interpretation of
vertex operator algebras}, Ph.D. thesis, Rutgers University, 1990.

\bibitem[H2]{H2} Y.-Z. Huang, {\it Geometric interpretation of vertex
operator algebras}, Proc. Natl. Acad. Sci. USA 88 (1991), 9964--9968.

\bibitem[H3]{H3} Y.-Z. Huang, {\it Two-Dimensional Conformal Geometry
and Vertex Operator Algebras}, Progress in Math., Vol. 148,
Birkh\"auser, Boston, 1997.

\bibitem[R]{R} C. Reutenauer, {\it Free Lie Algebras}, London
Math. Soc.  Monographs, New Series, 7, Clarendon Press, Oxford, 1993.

\bibitem[V]{V} V.S. Varadarajan, {\it Lie Groups, Lie Algebras, and
Their Representations}, Grad. Texts in Math., Vol. 120,
Springer-Verlag, New York, 1974.

\end{thebibliography}
\end{document}